\documentclass[12pt]{article}
\usepackage{latexsym,amsmath,amssymb}
\usepackage{graphicx}
\begin{document}
\title{{\large Every 4-regular graph is acyclically edge-6-colorable}}

\author{Weifan Wang\thanks{Research supported partially by
NSFC(No.11071223) and ZJNSF(No.Z6090150); Corresponding author.
Email: wwf@zjnu.cn.}, \ Qiaojun Shu\\
\bigskip
 \normalsize Department of Mathematics,  Zhejiang Normal University, Jinhua 321004, China\\
 Yiqiao Wang\\
 \normalsize School of   Management, Beijing University of Chinese
 Medicine,   Beijing 100029, China}

\maketitle
\newtheorem{theorem}{Theorem}
\newtheorem{lemma}{Lemma}
\newtheorem{claim}{Claim}
\newtheorem{corollary}[theorem]{Corollary}
\newtheorem{proposition}[theorem]{Proposition}
\newtheorem{conjecture}{Conjecture}
\newcommand{\qed}{\hfill $\Box$ }
\newcommand{\proof}{\noindent{\bf Proof.}\ \ }
\baselineskip=18pt
\parindent=0.5cm

\begin{abstract}

\baselineskip=16pt

An acyclic edge coloring of a graph $G$ is a proper edge coloring
such that no bichromatic cycles are produced. The acyclic chromatic
index $a'(G)$ of $G$ is the smallest integer $k$ such that $G$ has
an acyclic edge coloring using $k$ colors. Fiam${\rm \check{c}}$ik
(1978) and later Alon, Sudakov and Zaks (2001) conjectured that
$a'(G)\le \Delta + 2$ for any simple graph $G$ with maximum degree
$\Delta$. Basavaraju and  Chandran (2009) showed that every graph
$G$ with $\Delta=4$, which is not 4-regular,
  satisfies the conjecture. In this paper, we settle the 4-regular
  case, i.e., we show that every 4-regular graph $G$ has $a'(G)\le
  6$.

\medskip


\end{abstract}

\section{Introduction}
Only simple graphs are considered in this paper. Let $G$ be a graph
with vertex set $V(G)$ and edge set $E(G)$. A {\em proper
edge-$k$-coloring} is a mapping $c: E(G) \to \{1,2,\ldots ,k\}$ such
that any two adjacent edges receive different colors. The graph $G$
is {\em edge-$k$-colorable} if it has an edge-$k$-coloring. The {\em
chromatic index} $\chi'(G)$ of $G$ is the smallest integer $k$ such
that $G$ is edge-$k$-colorable. A proper edge-$k$-coloring   of $G$
is called {\em acyclic} if there are no bichromatic cycles in $G$,
i.e., the union of any two color classes induces a subgraph of $G$
that is a forest. The {\em acyclic chromatic index} of $G$, denoted
$a'(G)$, is the smallest integer $k$ such that $G$ is acyclically
edge-$k$-colorable.

Let $\Delta(G)$ ($\Delta$ for short) denote the maximum degree of a
graph $G$. By Vizing's Theorem \cite{Vizing64}, $\Delta \le
\chi'(G)\le \Delta +1$. Thus, it holds trivially  that $a'(G)\ge
\chi'(G)\ge \Delta$. Fiam${\rm \check{c}}$ik \cite{F78} and later
Alon, Sudakov and Zaks \cite{ABZ01} made the following conjecture:

\begin{conjecture}\label{con1}
For any graph $G$, $a'(G)\le \Delta+2$.
\end{conjecture}

Using probabilistic method, Alon, McDiarmid and Reed \cite{AMR91}
proved that $a'(G)\le 64\Delta$ for any graph $G$. This bound has
been recently improved to that  $a'(G)\le 16\Delta$ in \cite{MR98},
that $a'(G)\le \lceil 9.62(\Delta-1)\rceil$ in \cite{nd}, and that
$a'(G)\le  4\Delta $ in \cite{es}.  In  \cite{NNCW05}, N$\check{\rm
e}$set$\check{\rm r}$il and Wormald proved that $a'(G)\le \Delta+1$
for a random $\Delta$-regular graph $G$. The acyclic edge coloring
of some special classes of graphs was also investigated, including
subcubic graphs \cite{ BLSC08,Sku04}, outerplanar graphs
\cite{HWLL10, MNCRS07}, series-parallel graphs \cite{HWLL09,ws},
2-degenerate graphs \cite{ML12}, and planar graphs \cite{BLSNFT11,
CHT09, DX10, FHN08, SW10, WSWW11, wsw, YHLLX09}. In particular,
Conjecture \ref{con1} was confirmed for planar graphs without
3-cycles \cite{SW11}, or without 4-cycles \cite {WSW10}, or without
5-cycles \cite{SWW11}. It was shown in \cite{BLSNFT11} that every
planar graph $G$ has $a'(G)\le \Delta+12$. This upper bound was
recently improved to that
 $a'(G)\le \Delta+7$  \cite{wsw}.

Basavaraju and  Chandran \cite{MLSC09} showed that if $G$ is a graph
with $\Delta=4$ and $|E(G)|\leq 2|V(G)|- 1$, then $a'(G)\leq 6$.
Equivalently, every non-regular graph $G$ with $\Delta=4$ satisfies
Conjecture 1. In this paper, we will settle the case where $G$ is
4-regular. Our result and the result of \cite{MLSC09} confirm
Conjecture 1   for graphs with $\Delta=4$.

\section{Acyclic chromatic indices}

In this section, we discuss the acyclic chromatic indices of
$4$-regular graphs. Before establishing our main result, we need to
introduce some notation.

Assume that $c$ is a partial acyclic edge-$k$-coloring of a graph
$G$ using the color set $C = \{1, 2, \ldots, k\}$. For a vertex
$v\in V(G)$, we use $C(v)$ to denote the set of colors assigned to
edges incident to $v$ under $c$. If the edges of a cycle $ux\ldots
vu$ are alternatively colored with colors $i$ and $j$, then we call
such cycle an {\em $(i, j)_{(u, v)}$-cycle}. If the edges of a path
$ux\ldots v$ are alternatively colored with colors $i$ and $j$, then
we call such path an {\em $(i, j)_{(u, v)}$-path}. In the following
Figures, black spots are pairwise distinct, whereas the others may
be not.

For simplicity, we use {\em $\{e_{1}, e_2, \cdots, e_m\} \to a$} to
express that all edges $e_{1}$, $e_2$, $\cdots$, $e_m$ are colored
or recolored with the color $a$. In particular, when $m=1$, we write
simply $e_1 \to a$. Moreover,  we use  $(e_1, e_2, \cdots, e_m)_c=
(a_1, a_2, \cdots, a_m)$ to denote the fact that $c(e_i)= a_i$ for
$i=1,2,\ldots, m$. Let  $(e_1, e_2, \cdots, e_n) \to (b_1, b_2,
\cdots, b_n)$ denote the fact that $e_i$ is colored or recolored
with the color $b_i$ for $i=1,2,\ldots,n$.

Now we give two useful lemmas, which will be used in the following
discussion.

%
%

\begin{lemma}\label{12}
Suppose that a   graph $G$ has an edge-6-coloring $c$ such that
$c(uv)=i$,  where $uv$ is an arbitrary edge of $G$. If $uv$ can be
recolored properly with the color $j \ (\ne i)$, then $G$ contains
no $(i,j)_{(u, v)}$-path under $c$.
\end{lemma}

Lemma \ref{12} holds obviously by the definition.

\begin{lemma}\label{exchange}
Suppose that a graph $G$ has an edge-6-coloring $c$. Let  $P=
uv_1v_2\ldots v_kv_{k + 1}$ be  an $(i, j)_{(u, v_{k + 1})}$-path in
$G$ with $c(uv_1)= i$ and $j\not\in C(u)$. If $w\not\in V(P)$, then
there does not exist an $(i,j)_{(u, w)}$-path in $G$ under $c$.

\end{lemma}

\proof  Otherwise, assume that there is an $(i, j)_{(u, w)}$-path
$Q=uw_1w_2\ldots w_m$ in $G$, where $w=w_m$ and $m\ge 1$. Note that
some $w_s$ may be identical to some $v_t$, $1\le s\le m-1$, $1\le
t\le k+1$. Since $j\notin C(u)$ and $w\not\in V(P)$, it is easy to
see that there exists  a vertex $v_t$ and some vertex $w_s$ such
that $v_tv_{t+1}$ and $v_tw_s$ have the same color $i$ or $j$, which
contradicts the fact that $c$ is  a proper edge coloring. \qed

\begin{theorem}\label{triangle-free}
If $G$ is a $4$-regular graph without 3-cycles,
 then $a'(G)\leq 6$.
\end{theorem}
\proof  By Vizing's theorem \cite{Vizing64}, $G$ admits a proper
edge-6-coloring $c$ using the color set $C=\{1,2,\ldots,6\}$. Let
$\tau(c)$ denote the number of bichromatic cycles in $G$ with
respect to the coloring $c$. If $\tau(c)=0$, then $c$ is an acyclic
edge-$6$-coloring of $G$, hence we are done. Otherwise, $\tau(c)>0$.
Let $B$ be a bichromaic cycle of $G$. We are going to show that
there is a proper  edge-6-coloring $c'$, formed by
 recoloring suitably some edges of $G$, such that
 $B$ is no longer a bichromatic cycle
 and on new bichromatic cycles are produced. Namely,
 $\tau(c')<\tau(c)$. Repeating this process, we finally obtain an
 acyclic edge-6-coloring $c^*$ of $G$.



To arrive at our conclusion, assume, w.l.o.g., that $B$ is a
$(4,1)_{(u,v)}$-cycle with $c(uv)=1$.  Let $u_1,u_2,u_3$ be the
neighbors of $u$ different from $v$, and $v_1,v_2,v_3$  the
neighbors of $v$ different from $u$. Since $G$ contains no
triangles, $u_1,u_2,u_3,v_1,v_2,v_3$ are pairwise distinct. Clearly,
we may assume that $c(uu_2)=c(vv_2)=4$. Since
  $2\leq |C(u)\cap C(v)|\leq 4$, the proof is split into the
  following three cases, as shown in Fig.\,1.

\medskip
\begin{figure}[hh]
\begin{center}

\includegraphics[width= 6 in]{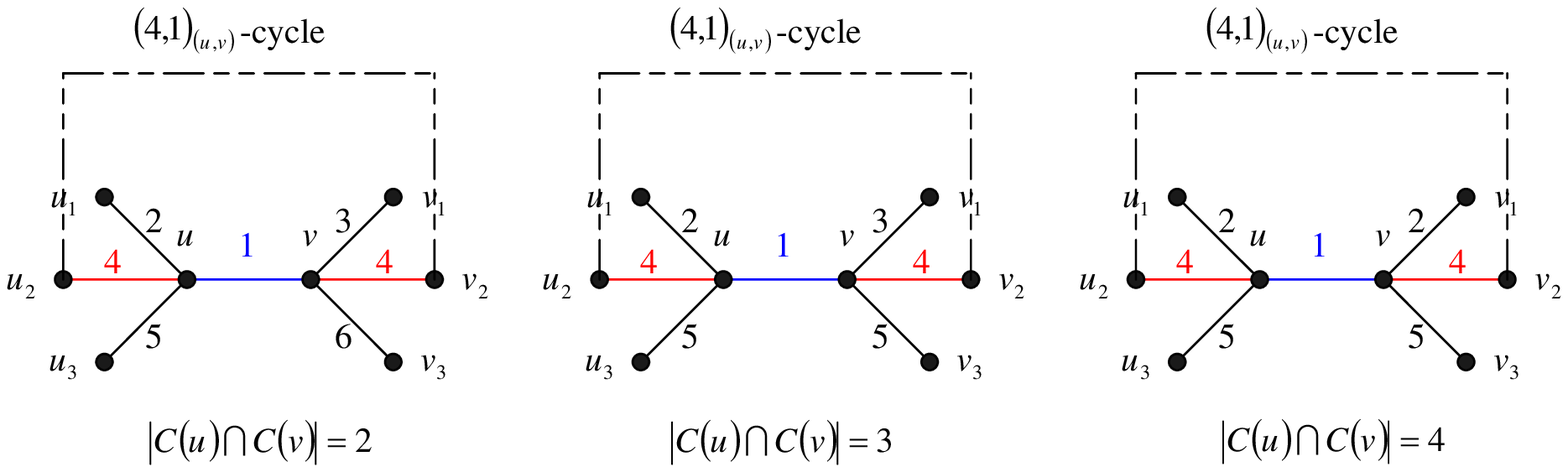}\\

\bigskip

{\rm Fig.\,1: Three cases on $|C(u)\cap C(v)|=2,3,4.$}
\end{center}
\end{figure}

 \noindent {\bf Case 1}\  $|C(u)\cap C(v)|= 2$, say $(uu_1, uu_3)_c=(2, 5)$ and $(vv_1,
vv_3)_c=(3, 6)$.

\medskip

By symmetry, we may assume that $C(u_2)\in
\{\{1,4,3,6\},\{1,4,2,6\},\{1,4,2,5\}\}$.

\medskip

 \noindent{\bf Case 1.1}\  $C(u_2)= \{1, 4, 3, 6\}$.
\medskip

Assume that $3\not\in C(u_1)\cup C(u_3)$. If $G$ contains no $(4,
3)_{(u_1, u_2)}$-path, then let $(uu_1, uv)\to (3, 2)$. Otherwise,
$G$ contains a $(4, 3)_{(u_1, u_2)}$-path and let $(uu_3, uv)\to (3,
5)$. Otherwise, assume that $3\in C(u_1)\cup C(u_3)$ and $6\in
C(u_1)\cup C(u_3)$. If $1\not\in C(u_1)$ and $G$ contains no $(1,
5)_{(u_1, u_3)}$-path, then let $(uu_1, uv)\to (1, 2)$. If $4\not\in
C(u_1)$ and $G$ contains no $(4, 5)_{(u_1, u_3)}$-path, then let
$(uu_1, uu_2)\to (4, 2)$. By Lemma \ref{exchange}, no new
bichromatic cycles are produced. Otherwise, assume that:

\medskip
\noindent ($*_{\textbf{1.1}}$)\

{\em
{\rm (a)}\ $\{1, 4\}\subseteq C(u_1)$, or
$G$ contains a $(5, i)_{(u_1, u_3)}$-path for any $i\in \{1,4\}\backslash C(u_1)$;

{\rm (b)}\ $\{1, 4\}\subseteq C(u_3)$, or $G$ contains a $(2, i)_{(u_1,
u_3)}$-path for any $i\in \{1, 4\}\backslash C(u_3)$;

{\rm (c)}\ $\{1, 4, 3, 6\}\subseteq C(u_1)\cup C(u_3)$;

{\rm (d)}\ $5\in C(u_1)$ if $\{1, 4\}\backslash C(u_1)\neq \emptyset$ and
$2\in C(u_3)$ if $\{1, 4\}\backslash C(u_3)\neq \emptyset$.}

\medskip
 \noindent{\bf (1.1.1)}\ Assume that $2\not\in C(u_3)$.
\medskip

By ($*_{\textbf{1.1}}$), $\{1, 4\}\subseteq C(u_3)$. If $1\not\in
C(u_1)$, then let $(uu_1, uu_3, uv)\to (1, 2, 5)$. If $4\not\in
C(u_1)$, then let $(uu_1, uu_2, uu_3)\to (4, 5, 2)$. Thus, $\{1,
4\}\subseteq C(u_1)$ and assume that $C(u_1)= \{1, 2, 3, 4\}$ and
$C(u_3) = \{5, 4, 1, 6\}$ since $\{3, 6\}\subseteq C(u_1)\cup
C(u_3)$ by ($*_{\textbf{1.1}}$). If $1\not\in C(v_1)$, then let
$(uu_2, uu_3)\to (5, 3)$. If $1\not\in C(v_3)$, then let $(uu_1,
uu_2)\to (6, 2)$. Otherwise, $1\in C(v_1)\cap C(v_3)$.

Suppose that we can recolor some edges in  $\{vv_1,vv_2,vv_3\}$ such
that $C(u)\cap C(v)\in \{\{1, 2\}, \{1, 5\}\}$, or $C(u)\cap C(v)\in
\{\{1, 2, 4\}, \{1, 5, 4\}\}$ and $c(vv_2)\neq 4$, and no new
bichromatic cycles are produced in $G-uv$. By Lemma \ref{exchange},
$B$ does not exist even if $4\in C(v)$. If $G$ contains a $(1,
i)_{(u, v)}$-cycle for some $i\in \{2, 5\}$, then let $(uu_1,
uu_3)\to (5, 2)$ and no new bichromatic cycles are produced. Next,
we claim that:

\medskip
\noindent ($\divideontimes_{\textbf{1}}$)\  \
{\em If we can
recolor some edges in  $\{vv_1,vv_2,vv_3\}$ such
 that $C(u)\cap
C(v)\in \{\{1, 2\}, \{1, 5\}\}$, or $C(u)\cap C(v)\in \{\{1, 2, 4\},
\{1, 5, 4\}\}$ and $c(vv_2)\neq 4$, and no new bichromatic cycles
are produced except the $(1, i)_{(u, v)}$-cycle for some $i\in \{2,
5\}$, we are done.}

$\bullet$ \ Assume that $C(v_2)= \{4, 1, 2, 5\}$.

If $G$ contains no $(3, i)_{(v_1, v_3)}$-path for some $i\in \{2, 5,
4\}\backslash C(v_3)$, then let $(vv_3, vv_2)\to (i, 6)$ by Lemma \ref{exchange}. Otherwise, $G$ contains a $(3,
i)_{(v_1, v_3)}$-path for every $i\in \{2, 5, 4\}\backslash C(v_3)$
and a $(6, i)_{(v_1, v_3)}$-path for every $i\in \{2, 5, 4\}\backslash
C(v_1)$ similarly. Thus, $3\in C(v_3)$ and $6\in C(v_1)$. However,
$\{2, 5, 4\}\subseteq C(v_1)\cup C(v_3)$, a contradiction.

$\bullet$ \ Assume that  $\{2, 5\}\backslash C(v_2)\neq \emptyset$
and $2\not\in C(v_2)$. (If $5\not\in C(v_2)$, we can first let
$(uu_1, uu_3)_c= (5, 2)$ and have a similar discussion.)

By symmetry, $C(v_2)= \{4, 1, 3, 5\}$ or $C(v_2)= \{4, 1, 3, 6\}$.

(i)\   $C(v_2)= \{4, 1, 3, 6\}$.

If $G$ contains neither a $(2, 3)_{(v_1, v_2)}$-path nor a $(2,
6)_{(v_2, v_3)}$-path, then let $vv_2\to 2$   by
($\divideontimes_{\textbf{1}}$). Otherwise, $G$ contains a $(2,
3)_{(v_1, v_2)}$-path and $2\in C(v_1)$. Assume that $2\not\in
C(v_3)$. Then first let $vv_3\to 2$. Next, if $G$ contains no $(4,
6)_{(u_2, v)}$-path, then let $uv\to 6$; otherwise, $G$ contains a
$(4, 6)_{(u_2, v)}$-path, we let $(uu_1, uu_3, uv)\to (6, 2,
5)$. Assume that $2\in C(v_3)$ and $5\in C(v_1)\cup
C(v_3)$ similarly (It follows that $\{2, 5\}\subseteq C(v_1)\cup C(v_3)$
if $G$ contains a $(2, 6)_{(v_2, v_3)}$-path). Then, let $(vv_1,
vv_2, vv_3)\to (4, 2, 3)$   by
($\divideontimes_{\textbf{1}}$).

(ii)\ $C(v_2)= \{4, 1, 3, 5\}$.

If $G$ contains no $(2, 3)_{(v_1, v_2)}$-path, then let $vv_2\to 2$.
Otherwise, $G$ contains a $(2, 3)_{(v_1, v_2)}$-path and $2\in
C(v_1)$. If $2\not\in C(v_3)$, then let $(vv_3, uv)\to (2, 6)$.
Otherwise, $2\in C(v_3)$. If $4\not\in C(v_1)$ and $G$ contains no
$(4, 6)_{(v_1, v_3)}$-path, then let $(vv_1, vv_2, uv)\to (4, 2,
3)$. Otherwise, $4\in C(v_1)$ or $G$ contains a $(4, 6)_{(v_1,
v_3)}$-path. If $4\not\in C(v_1)$ and $G$ contains a $(4, 6)_{(v_1,
v_3)}$-path, i.e., $C(v_1)= \{3, 1, 2, 6\}$, $C(v_3)= \{6, 1, 2,
4\}$, then let $(vv_1, vv_2, vv_3, uv)\to (4, 6, 5, 3)$. Otherwise,
$C(v_1)= \{3, 1, 2, 4\}$. If $5\not\in C(v_3)$, then let $(vv_2,
vv_3)\to (6, 5)$; otherwise, $C(v_3)= \{6, 1, 2, 5\}$,  let
$(vv_2, vv_3)\to (6, 4)$ and  we are done by
($\divideontimes_{\textbf{1}}$) or Lemma \ref{exchange}.

\medskip
 \noindent{\bf (1.1.2)}\ Assume that $2\in C(u_3)$ and $5\in C(u_1)$.
\medskip

Since $\{1, 4, 3, 6\}\subseteq C(u_1)\cup C(u_3)$ by
($*_{\textbf{1.1}}$), we conclude that $C(u_1)= \{2, 5, 1, 3,
5\}$ and $C(u_3)= \{5, 2, 4, 6\}$, or $C(u_1)= \{2, 5, 1, 4\}$ and
$C(u_3)= \{5, 2, 3, 6\}$. Let $x_1,x_2,x_3\ne u$ be the other three
neighbors of $u_1$  and $c(u_1x_1)= 1$. Note that
$G$ contains a $(1, 2)_{(x_1, u_3)}$-path and $2\in C(x_1)$ by
($*_{\textbf{1.1}}$).

$\bullet$ \ Assume that $C(u_1)= \{2, 5, 1, 3\}$, $C(u_3)= \{5, 2,
4, 6\}$, and $c(u_1x_2)= 3$

If $G$ contains no $(5, 6)_{(u_1, u_3)}$ -path, then let $(uu_1,
uv)\to (6, 2)$. Otherwise, $G$ contains a $(5, 6)_{(u_1,
u_3)}$-path. If $G$ contains no $(3, 2)_{(x_2, u_3)}$-path, then let
$(uu_2, uu_3)\to (5, 3)$; otherwise, assume that $G$ contains a $(3,
2)_{(x_2, u_3)}$ and $2\in C(x_2)$. If $\{4, 6\}\backslash
C(x_i)\neq \emptyset$ for some $i\in \{1, 2\}$, then let $(uu_2,
uu_3, uv)\to (5, c(u_1x_i), 2)$ and recolor $uu_1$ with a color in
$\{4, 6\}\backslash C(x_i)$. Otherwise, $C(x_1)= \{1, 2, 4, 6\}$ and
$C(x_2)= \{3, 2, 4, 6\}$. Then, we switch the colors of $u_1x_1$ and
$u_1x_2$, and let $(uu_3, uv)\to (1, 5)$.

$\bullet$ \ Assume that $C(u_1)= \{2, 1, 5, 4\}$, $C(u_3)= \{5, 2,
3, 6\}$, and $c(u_1x_2)= 4$.

Note that $G$ contains a $(4, 2)_{(u_1, u_3)}$-path by
($*_{\textbf{1.1}}$) and then $2\in C(x_2)$. Assume that $3\not\in
C(x_2)$. If $G$ contains no $(3, 5)_{(u_1, u_3)}$-path, then let
$(uu_1, uv)\to (3, 2)$; otherwise, $G$ contains a $(3, 5)_{(u_1,
u_3)}$-path and then let $(uu_1, uu_2, uu_3, uv)\to (3, 5, 4, 2)$.
Otherwise, $3\in C(x_2)$ and $6\in C(x_2)$ similarly. Thus, $C(x_2)=
\{4, 2, 3, 6\}$. Assume that $3\not\in C(x_1)$ and first let
$uu_1\to 3$. Next, if $G$ contains no $(3, 4)_{(u_1, u_2)}$-path,
then let $(uu_3, uv)\to (1, 2)$; otherwise, $G$ contains a $(3,
4)_{(u_1, u_2)}$-path, we let $(uu_2, uu_3, uv)\to (2, 4, 1)$.
Otherwise, $3\in C(x_1)$ and $6\in C(x_1)$ similarly. Thus, $C(x_1)=
\{1, 2, 3, 6\}$. Then, switch the colors of $u_1x_1$ and $u_1x_2$
and let $(uu_3, uv)\to (1, 5)$.

By ($*_{\textbf{1.1}}$) and Lemma \ref{exchange}, no new bichromatic
cycles are produced, as shown in Fig.\,2.

\medskip
\begin{figure}[hh]
\begin{center}

\includegraphics[width= 5 in]{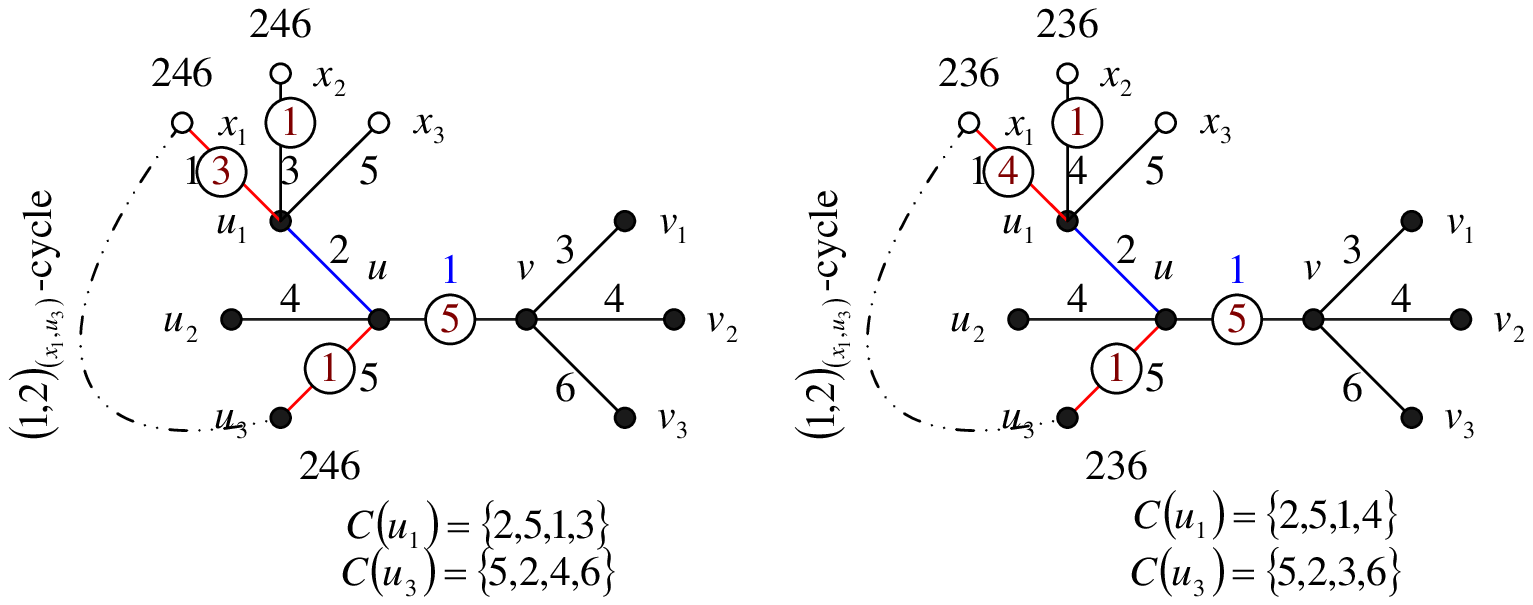}\\

\bigskip

{\rm Fig.\,2: The coloring in Case 1.1.2.}
\end{center}
\end{figure}

\medskip
 \noindent{\bf Case 1.2}\  $C(u_2)= \{1, 4, 2, 6\}$.
\medskip

If $C(v_2)= \{1, 4, 2, 5\}$,  we reduce the proof to  Case 1.1.
Thus, $\{2, 5\}\backslash C(v_2)\neq \emptyset$. If $G$ contains no
$(i, 2)_{(u_1, u_3)}$-path for some $i\in \{1, 3\}\backslash
C(u_3)$, then let $(uu_3, uv)\to (i, 5)$. Otherwise, assume that:

\medskip
\noindent ($*_{\textbf{1.2}}$)\

{\em
{\rm (a)}\ $\{1, 3\}\subseteq C(u_3)$, or
$G$ contains a $(i, 2)_{(u_1, u_3)}$-path for any $i\in \{1,
3\}\backslash C(u_3)$;

{\rm (b)}\ $\{1, 3\}\subseteq C(u_3)\cup C(u_1)$, and $2\in C(u_3)$ if
$\{1, 3\}\backslash C(u_3)\neq\emptyset$.}

\medskip
\noindent{\bf Case 1.2.1}\  $2\not\in C(u_3)$.
\medskip

Then $\{1, 3\}\subseteq C(u_3)$ by ($*_{\textbf{1.2}}$). If
$4\not\in C(u_1)$, then let $(uu_1, uu_2, uu_3)\to (4, 5, 2)$. If
$1\not\in C(u_1)$ and $G$ contains no $(4, 2)_{(u_2, u_3)}$-path,
then let $(uu_1, uu_3, uv)\to (1, 2, 5)$. Otherwise, $1\in C(u_1)$
or $G$ contains a $(4, 2)_{(u_2, u_3)}$-path.

$\bullet$ \ Assume that $C(u_3)= \{5, 1, 3, 6\}$.

Clearly, $1\in C(u_1)$ and first let $uu_3\to 4$. Next, if $5\not\in
C(u_1)$, then let $uu_2\to 5$ and we are done by Lemma
\ref{exchange}. Otherwise, $C(u_1)= \{2, 1, 4, 5\}$, we let $(uu_2,
uv)\to (3, 5)$.

$\bullet$ \ Assume that $C(u_3)=\{5, 1, 3, 4\}$.

If $G$ contains no $(4, 6)_{(u_2, u_3)}$-path, then let $(uu_3, uv)
\to (6, 5)$. Otherwise, assume that $G$ contains a $(4, 6)_{(u_2,
u_3)}$-path.  We have to handle two possibilities:

(i)\   $1\not\in C(u_1)$.

It follows that $G$ contains a $(4, 2)_{(u_2, u_3)}$-path. If
$5\not\in C(u_1)$, then let $(uu_1, uv)\to (1, 2)$. If $6\not\in
C(u_1)$, then let $(uu_1, uv)\to (6, 2)$. Otherwise, $C(u_1)= \{2,
4, 5, 6\}$ and first let $uu_3\to 6$. If $G$ contains no $(3,
6)_{(u_1, u_3)}$-path, let $(uu_1, uu_2, uv)\to (3, 5, 2)$;
otherwise, let $(uu_2,
uv)\to (3, 5)$.

(ii)\   $1\in C(u_1)$.

First assume that $3, 5\not\in C(u_1)$. Let  $uu_1\to 3$.
If $G$ contains no $(2, 4)_{(u_2, v_2)}$-path, we let $uv\to 2$; otherwise,
  let $(uu_3,
uv)\to (2, 5)$.

Next assume that $\{3, 5\}\cap C(u_1)\neq
\emptyset$ and $6\not\in C(u_1)$.  If $3\not\in C(u_1)$, then let
$(uu_1, uu_2, uu_3, uv)\to (3, 5, 6, 2)$.  So assume that $3\in C(u_1)$, i.e.,
  $C(u_1)= \{2, 1, 4, 3\}$.  Let $uu_1\to 6$. If $G$ contains no $(4,
2)_{(u_2, v_2)}$-path, further let $uv\to 2$; otherwise,  $(uu_3, uv)\to (2, 5)$.

\medskip
\noindent{\bf Case 1.2.2}\   $2\in C(u_3)$ and $G$ contains a $(2,
i)_{(u_1, u_3)}$-path for any $i\in \{1, 3\}\backslash C(u_3)$.
\medskip

{\bf (1.2.2.1)}  \ Assume that $G$ contains no $(2, 3)_{(u_1, u_2)}$-path.
\medskip

If $G$ contains no $(1, 3)_{(u_2, v_1)}$-path, then let $uu_2\to 3$.
Otherwise, $G$ contains a $(1, 3)_{(u_2, v_1)}$-path and $1\in
C(v_1)$.  If $C(v_1)= \{3, 1, 2, 5\}$, then by Case 1.1, we can
break that $(1, 3)_{(u, v)}$-cycle by setting $uu_2\to 3$ and hence break
the bichromatic cycle $B$ such that no new bichromatic cycles are produced.
Thus, $\{2, 5\}\backslash C(v_1)\neq \emptyset$.

Assume that there exists $a\in \{2, 5\}\backslash (C(v_1)\cup
C(v_2))$. If $G$ contains no $(6, a)_{(v_1, v_3)}$-path, then let
$(vv_1, uv)\to (a, 3)$; otherwise,  let $(vv_2, uu_2, uv)\to (a, 3,
4)$. Thus, assume that $\{2, 5\}\subseteq C(v_1)\cup C(v_2)$ and
$2\in C(v_1)\backslash C(v_2)$, $5\in C(v_2)\backslash C(v_1)$ (if
$5\in C(v_1)\backslash C(v_2)$ and $2\in C(v_2)\backslash C(v_1)$,
we can
 let $uu_2\to 3$ and  have a similar discussion).

If $G$ contains neither a $(2, 3)_{(v_1, v_2)}$-path nor a $(2,
6)_{(v_2, v_3)}$-path, then let $(vv_2$, $uu_2$, $uv)$ $\to$ $(2, 3,
4)$. Otherwise, $G$ contains either a $(2, 3)_{(v_1, v_2)}$-path or
a $(2, 6)_{(v_2, v_3)}$-path.

$\bullet$ \ Assume that $G$ contains a $(2, 3)_{(v_1, v_2)}$-path, implying
$C(v_2)= \{4, 1, 5, 3\}$.

If $\{1, 2\}\backslash$ $C(v_3)\neq\emptyset$, then let $uv\to 6$
and $vv_3\to \alpha\in \{1, 2\}\backslash C(v_3)$; otherwise, $\{1,
2\}\subseteq C(v_3)$. If $4\not\in C(v_1)$ and $G$ contains no $(4,
6)_{(v_1, v_3)}$-path, then let $(vv_1, vv_2, uv)\to (4, 2, 3)$;
otherwise, $4\in C(v_1)$, or $G$ contains a $(4, 6)_{(v_1,
v_3)}$-path and $C(v_1)= \{3, 1, 2, 6\}$, $C(v_3)= \{6, 1, 2, 4\}$.

If $4\not\in C(v_1)$, then let $(vv_1, uv)\to (5, 3)$. Otherwise,
$C(v_1)= \{3, 1, 2, 4\}$. If $G$ contains no $(4, 5)_{(v_1,
v_2)}$-path, then let $(vv_1, uv)\to (5, 3)$; otherwise, $G$
contains a $(4, 5)_{(v_1, v_2)}$-path. If $5\not\in C(v_3)$, then
let $(vv_3, uv)\to (5, 6)$. Otherwise, $C(v_3)= \{6, 1, 2, 5\}$ and
let $\{uu_2, vv_3\}\to 3$, $(vv_1, vv_2, uv)\to (5, 6, 4)$.

$\bullet$ \  \ Assume that $G$ contains a $(2, 6)_{(v_2, v_3)}$-path, implying
$C(v_2)= \{4, 1, 5, 6\}$.

First, assume that $C(v_1)= \{3, 1, 2, 4\}$. If $G$ contains no $(4,
5)_{(v_1, v_2)}$-path,  let $(vv_1, uv)\to (5, 3)$; otherwise,
$G$ contains a $(4, 5)_{(v_1, v_2)}$-path. If $5\not\in C(v_3)$,
 let $(vv_3, uu_2, uv)\to (5, 3, 6)$. If $3\not\in C(v_3)$,
let $\{uu_2, vv_3\}\to 3$ and $(vv_2, vv_1, uv)\to (2, 6, 4)$.
Otherwise, $C(v_3)= \{6, 2, 5, 3\}$ and let $(vv_1, vv_2, vv_3,
uv)\to (6, 2, 4, 3)$.

Next, assume that $C(v_1)= \{3, 1, 2, 6\}$. If $3\not\in C(v_3)$,
then let $\{vv_3, uu_2\}\to 3$ and $(vv_1, vv_2, uv)\to (4, 2, 6)$.
If $4\not\in C(v_3)$, then let $(vv_2, vv_3, uv)\to (2, 4, 6)$.
Otherwise, $C(v_3)= \{6, 2, 3, 4\}$ and let $(vv_1, vv_3, uv)\to (5,
1, 3)$.

\medskip
{\bf (1.2.2.2)} \  Assume that $G$ contains a $(2, 3)_{(u_1, u_2)}$-path.
\medskip

Then $3\in C(u_3)$ and $1\in C(u_1)\cup C(u_3)$ by
($*_{\textbf{1.2}}$). If $4\not\in C(u_1)\cup C(u_3)$,  let
$(uu_1, uu_2, uv)\to (4, 3, 2)$. Otherwise, $4\in C(u_1)\cup
C(u_2)$. If $4\not\in C(u_3)$ and $6\not\in C(u_1)\cup C(u_3)$,
let $(uu_3, uv)\to (6, 5)$. Otherwise, $4\in C(u_3)$ or $6\in
C(u_1)\cup C(u_3)$.

$\bullet$ \  Assume that $C(u_3)= \{5, 2, 3, 1\}$ and $C(u_1)= \{2, 3, 6,
4\}$.

If $G$ contains no $(4, 2)_{(u_2, v_2)}$-path,   let $(uu_1, uv)
\to (1, 2)$; otherwise,   let $(uu_1, uu_2, uu_3,$ $ uv)\to (5, 3, 4, 2)$.

$\bullet$ \  Assume that $C(u_3)= \{5, 2, 3, 6\}$ and $C(u_1)= \{2, 3, 1,
4\}$.

If $G$ contains no $(4, 2)_{(u_2, v_2)}$-path,   let $(uu_1,
uu_3, uv) \to (5, 1, 2)$; otherwise,   let $(uu_1, uu_2, uu_3, uv)\to (5, 3, 4, 2)$.

$\bullet$ \  Assume that $C(u_3)= \{5, 2, 3, 4\}$ and $\{2, 3,
1\}\subseteq C(u_1)$.

If $4, 5\not\in C(u_1)$,  let $(uu_1, uu_2, uv)\to (4, 3, 2)$.
If $5, 6\not\in C(u_1)$,   let $(uu_3, uu_2)\to (6, 5)$.
Otherwise, $\{4, 5\}\cap C(u_1)\neq \emptyset$ and $\{5, 6\}\cap
C(u_1)\neq \emptyset$. Thus, $C(u_1)= \{2, 3, 1, 5\}$. If $G$
contains no $(5, 6)_{(u_1, v_3)}$-path, let $(uu_1, uu_3, uv)\to (6,
1, 5)$; otherwise,   let
$(uu_1, uu_2, uu_3, uv)\to (6, 5, 1, 2)$.

\medskip
 \noindent{\bf Case 1.3}\   $C(u_2)= \{1, 4, 2, 5\}$.
\medskip

If $C(v_2)\neq \{1, 4, 3, 6\}$, then we can reduce   to   Cases 1.1
or  1.2. Thus, $C(v_2)= \{1, 4, 3, 6\}$. If $\{1, 3, 6\}\backslash
C(u_1)\neq \emptyset$ and $G$ contains no $(5, i)_{(u_1, u_3)}$-path
for some $i\in \{1, 3, 6\}\backslash C(u_1)$, then let $(uu_1,
uv)\to (i, 2)$. Otherwise, assume that:

\medskip
\noindent ($*_{\textbf{1.3}}$)\

(a)\ $\{1, 3, 6\}\subseteq C(u_1)$, or $G$ contains a $(5, i)_{(u_1,
u_3)}$-path for any $i\in \{1, 3, 6\}\backslash C(u_1)$;

(b)\ $\{1, 3, 6\}\subseteq C(u_3)$, or $G$ contains a $(2, i)_{(u_1,
u_3)}$-path for any $i\in \{1, 3, 6\}\backslash C(u_3)$;

(c)\ $\{1, 2, 5\}\subseteq C(v_3)$, or $G$ contains a $(3, i)_{(v_1,
v_3)}$-path for any $i\in \{1, 2, 5\}\backslash C(v_3)$;

(d)\ $\{1, 2, 5\}\subseteq C(v_1)$ or $G$ contains a $(6, i)_{(v_1,
v_3)}$-path for any $i\in \{1, 2, 5\}\backslash C(v_1)$;

(e)\ $\{1, 3, 6\}\subseteq C(u_1)\cup C(u_3)$ and $\{1, 2,
5\}\subseteq C(v_1)\cup C(v_3)$;

(f)\ $5\in C(u_1)$ if $\{1, 3, 6\}\backslash C(u_1)\neq \emptyset$,
$2\in C(u_3)$ if $\{1, 3, 6\}\backslash C(u_3)\neq \emptyset$, $6\in
C(v_1)$ if $\{1, 2, 5\}\backslash C(v_1)\neq \emptyset$, and $3\in
C(v_3)$ if $\{1, 2, 5\}\backslash C(v_3)\neq \emptyset$.

\medskip
 \noindent{\bf Case 1.3.1}\   $5\not\in C(u_1)$.
\medskip

We see that $C(u_1)= \{2, 1, 3, 6\}$ by ($*_{\textbf{1.3}}$). Note
that $G$ contains a $(2, i)_{(u_1, u_3)}$-path for any $i\in \{3,
6\}\backslash C(u_3)$. If $\{1, 3, 6\}\backslash
C(u_3)\neq\emptyset$, then let $(uu_1, uv)\to (5, 2)$ and $uu_3\to
\beta\in \{1, 3, 6\}\backslash C(u_3)$. So assume that $C(u_3)= \{5,
1, 3, 6\}$.

First assume that $G$ contains no $(3, 2)_{(u_2, v_1)}$-path. If $G$
contains no $(3, 5)_{(u_2, u_3)}$-path, let $(uu_1, uu_2, uv) \to
(4, 3, 2)$; otherwise,  let $(uu_1, uu_2, uu_3, uv)\to (5, 3, 4,
2)$.

Next assume that $G$ contains a $(3, 2)_{(u_2, v_1)}$-path and a
$(6, 2)_{(u_2, v_3)}$-path. If $G$ contains no $(1, 2)_{(u_1,
v_2)}$-path,  let $vv_2\to 2$; otherwise, let $\{vv_2, uu_3\}\to 2$
and $uu_1\to 5$.

\medskip
 \noindent{\bf Case 1.3.2}\   $5\in C(u_1)$.
\medskip

Furthermore, we may assume that $2\in C(u_3)$, $3\in C(v_3)$, and
$6\in C(v_1)$. Since $\{1, 3, 6\}\subseteq C(u_1)\cup C(u_3)$ by
($*_{\textbf{1.3}}$), we suppose that $1\in C(u_1)$ and it suffices
to consider the following two subcases by symmetry:


$\bullet$ \ $C(u_1)= \{2, 5, 1, 3\}$, $6\in C(u_3)$ and $4\not\in
C(u_3)$.

By ($*_{\textbf{1.3}}$), $G$ contains a $(5, 6)_{(u_1, u_3)}$-path.
If $\{1, 2\}\backslash C(v_3)\neq \emptyset$,   let $(uu_1, uu_2)
\to (4, 6)$ and $uv\to \gamma\in \{1, 2\}\backslash C(v_3)$.
Otherwise, $C(v_3)= \{3, 6, 1, 2\}$ and $G$ contains a $(3,
5)_{(v_1, v_3)}$-path by ($*_{\textbf{1.3}}$), hence  let $(uu_2,
vv_2, uv)\to (6, 5, 4)$.

$\bullet$ \  $C(u_1)= \{2, 5, 1, 3\}$ and $C(u_3)= \{5, 2, 6, 4\}$;
or $C(u_1)= \{2, 5, 1, 4\}$ and $C(u_3)= \{5, 2, 3, 6\}$.

By ($*_{\textbf{1.3}}$), $G$ contains a $(5, 6)_{(u_1, u_3)}$-path,
and moreover $G$ contains a $(2, 3)_{(u_1, u_3)}$-path when $C(u_1)=
\{2, 5, 1, 3\}$, a $(5, 3)_{(u_1, u_3)}$-path when $C(u_1)= \{2, 5,
1, 4\}$. If $1\not\in C(v_i)$ for some $i\in \{1, 3\}$, then let
$uu_2\to c(vv_i)$. Otherwise, $1\in C(v_1)\cap C(v_3)$. Since $\{2,
5\}\subseteq C(v_1)\cup C(v_3)$ and noting that $1\not\in C(u_3)$,
we have to consider the following two possibilities:  If $C(v_1)=
\{3, 6, 1, 2\}$ and $C(v_3)= \{3, 6, 1, 5\}$,  let $(vv_2, vv_3)\to
(5, 4)$. If $C(v_1)= \{3, 6, 1, 5\}$ and $C(v_3)= \{3, 6, 1, 2\}$,
 let $(vv_1, vv_2)\to (4, 5)$.

\medskip

 \noindent {\bf Case 2}\  $|C(u)\cap C(v)|= 3$, say $(uu_1, uu_3)_c= (2, 5)$ and $(vv_1,
vv_3)_c=(3,5)$.

\medskip

If $C(u_2)= \{1, 4, 3, 6\}$ or $C(v_2)= \{1, 4, 2, 6\}$, then the
proof can be reduced to  Case 1. If $G$ contains neither a $(4,
6)_{(u_2, v)}$-path nor a $(5, 6)_{(u_2, v)}$-path, then let $uv\to
6$. Otherwise, assume that:

\medskip

\noindent ($*_{\textbf{2.1}}$)\ $G$ contains either a $(4, 6)_{(u,
v)}$-path or a $(5, 6)_{(u, v)}$-path.
\medskip

Suppose that we can recolor some edges in $G$ such that $|C(u)\cap
C(v)|= 2$, and no new bichromatic cycles are produced in $G-uv$.
Then, by the discussion of Case 1.1, if $G$ contains an $(i, j)_{(u,
v)}$-cycle for some $i, j\in C(u)\cap C(v)$, then we can break this
$(i, j)_{(u, v)}$-cycle and hence break the cycle $B$ such that no
new bichromatic cycles are produced. That is, we have the following
statement:

\medskip
\noindent ($\divideontimes_{\textbf{2}}$)\ If we can recolor some
edges in  $G$ such that $|C(u)\cap C(v)|= 2$, and no new bichromatic
cycles are produced in $G-uv$, then we are done.
\medskip

If $6\not\in C(u_3)$ and $G$ contains no $(2, 6)_{(u_1, u_3)}$-path,
 let $uu_3\to 6$; if $1\not\in C(u_3)$ and $G$ contains no $(1,
2)_{(u_1, u_3)}$-path,  let $(uu_3, uv)\to (1, 6)$. Clearly,
$|C(u)\cap C(v)|= 2$ and no new bichromatic cycles are produced in
$G-uv$. By ($\divideontimes_{\textbf{2}}$), we  complete the proof.
Similarly, if $6\not\in C(u_2)$ and $G$ contains no $(2, 6)_{(u_1,
u_2)}$-path, we can let $uu_2\to 6$. Hence, assume that:

\medskip
\noindent ($*_{\textbf{2.2}}$)\

(a)\ $6\in C(u_2)$, or $G$ contains a $(2, 6)_{(u_1, u_2)}$-path;

(b)\ $6\in C(v_2)$, or $G$ contains a $(3, 6)_{(v_1, v_2)}$-path;

(c)\  $\{1, 6\}\subseteq C(u_3)$, or $G$ contains a $(2, i)_{(u_1,
u_3)}$-path for each  $i\in \{1, 6\}\backslash C(u_3)$;

(d)\ $\{1, 6\}\subseteq C(v_3)$, or $G$ contains a $(3, i)_{(v_1,
v_3)}$-path for each  $i\in \{1, 6\}\backslash C(v_3)$;

(e)\ $\{1, 6\}\subseteq C(u_1)\cup C(u_3)$ and $\{1, 6\}\subseteq
C(v_1)\cup C(v_3)$;

(f)\ $2\in C(u_3)$ if $\{1, 6\}\backslash C(u_3)\neq \emptyset$ and
$3\in C(u_3)$ if $\{1, 6\}\backslash C(v_3)\neq \emptyset$.
\medskip

It follows that $\{2, 6\}\cap C(u_2)\neq \emptyset$ and $\{3,
6\}\cap C(v_2)\neq \emptyset$. Thus, $C(u_2) \backslash \{1, 4\}\in
\{\{2, 3\}, \{2, 5\}, \{2, 6\},\{5, 6\}\}$ and $C(v_2) \backslash
\{1, 4\}\in \{\{2, 3\}, \{3, 5\}, \{3, 6\}, \{5, 6\}\}$.

\medskip
 \noindent{\bf Case 2.1}\ $C(u_2)\backslash \{1, 4\}= \{2, 3\}$.
\medskip

By ($*_{\textbf{2.1}}$) and ($*_{\textbf{2.2}}$), $G$ contains a
$(5, 6)_{(u, v)}$-path with $6\in C(u_3)\cap C(v_3)$,  or  a
$(2, 6)_{(u_1, u_2)}$-path with  $6\in C(u_1)$   and $1\in C(u_3)$, or
 a $(1, 2)_{(u_1, u_3)}$-path with $2\in C(u_3)$ and $1\in
C(u_1)$. If $4\not\in C(u_1)$ and $G$ contains no $(4, 5)_{(u_1,
u_3)}$-path,   let $(uu_1, uu_2)\to (4, 6)$. Otherwise, $4\in
C(u_1)$, or $G$ contains a $(4, 5)_{(u_1, u_3)}$-path and $5\in
C(u_1)$, $4\in C(u_3)$. We need to consider the following
subcases.

\medskip
 \noindent{\bf Case 2.1.1} \   $4\in C(u_1)$ and $1\in C(u_3)$.
\medskip

$\bullet$ Assume that $C(u_1)= \{2, 6, 4, 3\}$.

If $2\not\in C(u_3)$,  let $(uu_1, uu_2, uu_3, uv)\to (1, 5, 2,
6)$; otherwise, $C(u_3)= \{5, 6, 1, 2\}$,  let $(uu_1, uu_2,
uu_3)\to (5, 6, 4)$.

$\bullet$ Assume that $3\not\in C(u_1)$.

If $2\not\in C(u_3)$ and $G$ contains no $(2, 4)_{(u_2, u_3)}$-path,
 let $(uu_3, uv)\to (2, 6)$ and $uu_1\to \alpha\in\{1,
5\}\backslash C(u_1)$. By ($*_{\textbf{2.1}}$) and Lemma
\ref{exchange}, no new bichromatic cycles are produced. Otherwise,
$2\in C(u_3)$ or $G$ contains a $(2, 4)_{(u_2, u_3)}$-path and $4\in
C(u_3)$.

(i)\ Assume that $C(u_3)= \{5, 1, 6, 4\}$ and $G$ contains a $(2,
4)_{(u_2, u_3)}$-path.

If $G$ contains no $(4, 3)_{(u_1, u_2)}$-path,   let $(uu_1,
uv)\to (3, 2)$; otherwise,   let $(uu_3, uv)\to (3, 2)$, $uu_1\to \beta\in\{1, 5\}\backslash
C(u_1)$.

(ii)\ Assume that $C(u_3)= \{5, 6, 1, 2\}$.

If $G$ contains no $(3, 6)_{(u_3, v_1)}$-path,  let $(uu_3,
uv)\to (3, 6)$. So assume that $G$ contains a $(3, 6)_{(u_3,
v_1)}$-path. It follows that $6\in C(v_2)$ by ($*_{\textbf{2.2}}$)
and the fact that $2\not\in C(v_2)$. If $C(u_1)= \{2, 6, 4, 1\}$,  let $(uu_1,
uu_2, uu_3, uv)\to (5, 6, 4, 2)$.  Otherwise, we derive that $C(u_1)= \{2, 6, 4,
5\}$. If $G$ contains no $(2, 3)_{(u_3, v_1)}$-path, then let
$(uu_1, uu_3, uv)\to (1, 3, 2)$. Thus, assume that $G$ contains a $(2,
3)_{(u_3, v_1)}$-path and $2\in C(v_1)$. If $G$ contains no $(2,
5)_{(v_2, v_3)}$-path, then let $vv_2\to 2$; otherwise, $G$ contains
a $(2, 5)_{(v_2, v_3)}$-path, which implies that $C(v_2)= \{4, 1, 6,
5\}$ and $2\in C(v_3)$. If $1, 5\not\in C(v_1)$,  let $(vv_1,
uv)\to (1, 3)$. Otherwise, $\{1, 5\}\cap C(v_1)\neq \emptyset$ and
then $4\not\in C(v_1)$ since $\{2, 6\}\subseteq C(v_1)$. If
$4\not\in C(v_3)$, then let $(vv_3, vv_2)\to (4, 2)$; otherwise, let
$(vv_1, vv_2, vv_3)\to (4, 2, 3)$.

\medskip
 \noindent{\bf Case 2.1.2} \   $4\notin C(u_1)$ or  $1\notin C(u_3)$.
\medskip

If $4\not\in C(u_1)$ and $1\in C(u_3)$, then $C(u_3)= \{5, 6, 1,
4\}$,  let $(uu_1, uu_2, uu_3)\to (4, 6, 2)$. If $4\in C(u_1)$ and
$1\not\in C(u_3)$, then $C(u_1)=\{2, 6, 1, 4\}$, let $(uu_1,
uu_3, uv)\to (5, 1, 6)$. Otherwise, assume that $4\not\in C(u_1)$
and $1\not\in C(u_3)$. Then $C(u_1)= \{2, 6, 1, 5\}$ and $C(u_3)=
\{5, 6, 2, 4\}$. If $G$ contains no $(3, 6)_{(u_1, v_1)}$-path,
let $(uu_1, uu_3, uv)\to (3, 1, 6)$. Otherwise, it follows that $6\in C(v_2)$ by
($*_{\textbf{2.2}}$) and $2\not\in C(v_2)$. Let $(uu_1, uu_3,
uv)\to (3, 1, 2)$.

\medskip
 \noindent{\bf Case 2.2}\ $C(u_2)\backslash \{1, 4\}= \{2, 5\}$.
\medskip

If $C(v_2)\backslash\{1, 4\}= \{2, 3\}$, then the proof is  reduced   to
 Case 2.1. Thus, $C(v_2)\backslash \{1, 4\}\neq\{2, 3\}$ and
$2\not\in C(v_2)$. By ($*_{\textbf{2.1}}$) and ($*_{\textbf{2.2}}$),
 $G$ contains a $(5, 6)_{(u_3,
v_3)}$-path with $6\in C(u_3)\cap C(v_3)$, or   $G$ contains a $(2, 6)_{(u_1,
u_2)}$-path with  $6\in C(u_1)$  and $1\in C(u_3)$, or a $(1, 2)_{(u_1,
u_3)}$-path with  $1\in C(u_1)$ and $2\in C(u_3)$.  If $4\not\in C(u_1)$
and $G$ contains no $(4, 5)_{(u_1, u_3)}$-path, then let $(uu_1,
uu_2)\to (4, 6)$. Otherwise,  $4\in C(u_1)$,  or $G$ contains no $(4,
5)_{(u_1, u_3)}$-path and $5\in C(u_1)$, $4\in C(u_3)$. We need to
consider the following two subcases.

\medskip
 \noindent{\bf Case 2.2.1}\ $4\in C(u_1)$ and $1\in C(u_3)$.
\medskip

$\bullet$ \ Assume that $C(u_1)= \{2, 6, 4, 1\}$.

If $2\not\in C(u_3)$,   let $(uu_1, uv)\to (3, 2)$. Otherwise,
$C(u_3)= \{5, 6, 1, 2\}$, let $(uu_1, uu_2, uu_3, uv)$ $\to (3, 6,
4, 2)$.

$\bullet$ \ Assume that $1\not\in C(u_1)$.

If $2, 4\not\in C(u_3)$,  let $(uu_1, uu_3, uv)\to (1, 2, 6)$.
Otherwise, $\{2, 4\}\cap C(u_3)\neq \emptyset$.

(i)\ Assume that $C(u_3)= \{5, 1, 6, 4\}$.

If $5\not\in C(u_1)$,  let $(uu_1, uv)\to (1, 2)$. Otherwise,
$C(u_1)= \{2, 6, 4, 5\}$, we let $(uu_1, uv)\to (3, 2)$.

(ii)\ Assume that $C(u_3)= \{5, 6, 1, 2\}$.

If $3\not\in C(u_1)$,  let $(uu_1, uu_2)\to (3, 6)$. Otherwise,
$C(u_1)= \{2, 6, 4, 3\}$.

If $G$ contains no $(2, 5)_{(u_3, v_3)}$-path,   let $(uu_1,
uv)\to (1, 2)$. So assume that $G$ contains a $(2, 5)_{(u_3, v_3)}$-path
and $2\in C(v_3)$. If $G$ contains no $(2, 3)_{(v_1, v_2)}$-path,
let $vv_2\to 2$. Thus, assume that $G$ contains a $(2, 3)_{(v_1,
v_2)}$-path and $2\in C(v_1)$, $3\in C(v_2)$. Note that $1\in
C(v_3)$ or $G$ contains a $(1, 3)_{(v_1, v_3)}$-path, and $\{1,
3\}\cap C(v_3)\neq \emptyset$ by ($*_{\textbf{2.2}}$) and then
$4\not\in C(v_3)$. If $4\not\in C(v_1)$, then let $(vv_1, vv_2)\to
(4, 2)$. Otherwise, $4\in C(v_1)$.

Assume that $1\not\in C(v_3)$ and then $C(v_3)= \{5, 6, 2, 3\}$,
$C(v_1)= \{3, 2, 4, 1\}$ by ($*_{\textbf{2.2}}$). If $6\not\in
C(v_2)$, let $vv_2\to 6$; otherwise, $C(v_2)= \{4, 1, 3, 6\}$,
let $(vv_1, vv_3, uv)\to (5, 1, 6)$. Next, assume that $C(v_3)= \{5,
6, 2, 1\}$. If $1\not\in C(v_1)$,   let $(vv_1, vv_3, uv)\to (1,
3, 6)$; otherwise, $C(v_1)= \{3, 2, 4, 1\}$. If $6\not\in C(v_2)$,
then let $vv_2\to 6$; otherwise, $C(v_2)= \{4, 1, 3, 6\}$,  let
$(vv_1, vv_3, uv)\to (5, 3, 6)$.

\medskip
 \noindent{\bf Case 2.2.2}\ $4\notin C(u_1)$ or  $1\notin C(u_3)$.
\medskip

If $4\not\in C(u_1)$ and $1\in C(u_3)$, then $C(u_3)= \{5, 6, 1,
4\}$,  let $(uu_1, uu_2, uu_3)\to (4, 6, 2)$. If $4\not\in C(u_1)$
and $1\not\in C(u_3)$, then $C(u_1)= \{2, 6, 1, 5\}$, $C(u_3)= \{5,
6, 2, 4\}$, and let $(uu_1, uu_3, uv)\to (3, 1, 2)$. Otherwise,
assume that $4\in C(u_1)$ and $1\not\in C(u_3)$. Then $C(u_1)=\{2,
6, 1, 4\}$. If $G$ contains no $(1, 3)_{(u_1, u_3)}$-path, then let
$(uu_1, uu_3, uv)\to (3, 1, 2)$; otherwise, let $(uu_1, uu_2, uu_3, uv)\to (5, 3, 1,
6)$.

\medskip
 \noindent{\bf Case 2.3}\ $C(u_2)\backslash \{1, 4\}= \{2, 6\}$.
\medskip

If $C(v_2)\backslash\{1, 4\}\in \{\{2, 3\}, \{3, 5\}\}$, then we reduce
the proof to  Cases 2.1 or  2.2. Thus, $C(v_2)\backslash
\{1, 4\}\in \{\{3, 6\}, \{5, 6\}\}$.

\medskip
 \noindent{\bf Case 2.3.1}\   $G$ contains no $(2, 3)_{(u_1,
u_2)}$-path.
\medskip

First let $uu_2\to 3$. If $G$ contains no $(1, 3)_{(u_2,
v_2)}$-path, then we are done. Otherwise, $G$ contains a $(1,
3)_{(u_2, v_2)}$-path and $|C(u)\cap C(v)|= 2$.

If $C(v_1)\backslash \{1, 3\}\not\in \{\{4, 6\}, \{5, 6\}\}$, then
by Case 1, Case 2.1 or Case 2.2, we can break this $(1, 3)_{(u,
v)}$-cycle and hence break the cycle  $B$ and no new bichromatic cycles
are produced. Thus, $C(v_1)\backslash \{1, 3\}\in \{\{4, 6\}, \{5,
6\}\}$.

Note that $G$ contains a $(3, i)_{(v_1, v_3)}$-path for any $i\in
\{1, 6\}\backslash C(v_3)$ by ($*_{\textbf{2.2}}$).  This implies that
$1\in C(v_3)$. If $6\not\in C(v_3)$,  let $uv\to 6$. Hence,
$\{1, 6\}\subseteq C(v_3)$.

$\bullet$ \ Assume that $C(v_2)= \{4, 1, 6, 3\}$. (If $C(v_1)= \{3, 1,
6, 4\}$, we have a similar discussion.)

If $4\not\in C(v_3)$,  let $(vv_2, uv)\to (2, 4)$. Otherwise,
$C(v_3)= \{5, 1, 6, 4\}$, we let $(vv_1, uu_2, uv)\to (2, 4, 3)$.

%

$\bullet$ \ Assume that $C(v_2)= \{4, 1, 6, 5\}$ and $C(v_1)= \{3,
1, 6, 5\}$.

If $3, 4\not\in C(v_3)$,  let $(vv_1, vv_2)= (4, 3)$. Otherwise,
$\{3, 4\}\cap C(v_3)\neq \emptyset$. If $C(v_3)= \{5, 1, 6, 3\}$,
 let $(vv_2, uv)= (2, 4)$; otherwise, $C(v_3)= \{5, 1, 6, 4\}$,
 let $(uu_2, vv_1, uv)= (4, 2, 3)$.

\medskip
 \noindent{\bf Case 2.3.2}\  $G$ contains a $(2, 3)_{(u_1,
u_2)}$-path and $3\in C(u_1)$.
\medskip

 {\bf (2.3.2.1)}\    $C(v_2)= \{4, 1, 5, 6\}$.

 $\bullet$ \  Assume that $G$ contains no $(2, 5)_{(v_2, v_3)}$-path. Let
$vv_2\to 2$. If $G$ contains no $(1, 2)_{(u, v)}$-path, then we are
done; otherwise, $G$ contains a $(1, 2)_{(u, v)}$-path,  we
can break this $(1, 2)_{(u, v)}$-cycle by Case 1 or Case 2.1
as $3\in C(u_1)$.

 $\bullet$ \  Assume that $G$ contains a $(2, 5)_{(v_2, v_3)}$-path. If $2\not\in
C(v_1)$ and $G$ contains no $(5, 3)_{(u_3, v_3)}$-path, then  let
$(vv_1, uv)\to (2, 3)$.  So assume that $2\in C(v_1)$ or $G$ contains a
$(5, 3)_{(u_3, v_3)}$-path. If $2\not\in C(v_1)$ and $G$ contains a
$(5, 3)_{(u_3, v_3)}$-path, let $(vv_1, vv_2)\to (2, 3)$. If
$4\not\in C(v_3)$ and $G$ contains no $(3, 4)_{(v_1, v_3)}$-path,
  let $(vv_3, vv_2)\to (4, 2)$. If $G$ contains no $(1, 2)_{(u,
v)}$-path,  we are done; otherwise, $G$ contains a $(1, 2)_{(u,
v)}$-path, we can break this $(1, 2)_{(u, v)}$-cycle by
Cases 1 or 2.1 as $3\in C(u_1)$. Otherwise, $2\in C(v_1)$,
and $4\in C(v_3)$ or $G$ contains a $(3, 4)_{(v_1, v_3)}$-path.
Together with ($*_{\textbf{2.2}}$), we have $3\in C(v_3)$, $\{1,
6\}\backslash C(v_3)\neq \emptyset$, and then $5\not\in C(v_1)$
since $\{1, 6, 4\}\subseteq C(v_1) \cup C(v_3)$. Let $(vv_1,
vv_2)\to (5, 3)$, $vv_3\to \beta\in \{1, 6\}\backslash C(v_1)$, and
$uv\to \{1, 6\}\backslash \beta$.  By
($\divideontimes_{\textbf{2}}$), we complete the proof.

 {\bf (2.3.2.2)}\   $C(v_2)= \{4, 1, 3, 6\}$.

Similarly, we can assume that $G$ contains a $(2, 3)_{(v_1,
v_2)}$-path and $2\in C(v_1)$.

First, assume that $2\not\in C(u_3)$,  and hence $\{1, 6\}\subseteq C(u_3)$ by
($*_{\textbf{2.2}}$).
 If $4\not\in C(u_1)\cup C(u_3)$,  let
$(uu_1, uu_2, uv)\to (4, 3, 2)$; otherwise, $4\in C(u_1)\cup
C(u_3)$. If $1\not\in C(u_1)$ and $G$ contains no $(2, 4)_{(u_2,
u_3)}$-path, then let $(uu_1, uu_3, uv)\to (1, 2, 6)$ and we are
done by ($\divideontimes_{\textbf{2}}$); otherwise, $1\in C(u_1)$ or
$G$ contains a $(2, 4)_{(u_2, u_3)}$-path. Assume that $C(u_3)= \{5,
1, 6, 3\}$ and then $\{1, 4\}\subseteq C(u_1)$ since $G$ contains no
$(2, 4)_{(u_2, u_3)}$-path and $4\in C(u_1)\cup C(u_3)$. Then, let
$(uu_1, uu_2, uu_3, uv)\to (5, 3, 4, 2)$. Otherwise, assume that
$C(u_3)=\{5, 1, 6, 4\}$. If $5\not\in C(u_1)$,  let $(uu_1,
uu_3, uv)\to (5, 3, 2)$; otherwise, $5\in C(u_1)$. Assume that
$1\not\in C(u_1)$. If $G$ contains no $(1, 5)_{(u_1, u_3)}$-path,
 let $(uu_1, uv)\to (1, 2)$; otherwise, $G$ contains a $(1,
5)_{(u_1, u_3)}$-path,  let $(uu_1, uu_2, uu_3, uv)\to (1, 5, 2,
6)$ and we are done by ($\divideontimes_{\textbf{2}}$). Otherwise,
$C(u_1)= \{2, 3, 5, 1\}$. If $G$ contains no $(5, 6)_{(u_1,
u_3)}$-path, then let $(uu_1, uv)\to (6, 2)$; otherwise, $G$
contains a $(5, 6)_{(u_1, u_3)}$-path, let $(uu_1, uu_2, uu_3)\to
(6, 5, 2)$ and we are done by ($\divideontimes_{\textbf{2}}$).

Next, assume that $2\in C(u_3)$ and $3\in C(v_3)$. Assume that
$4\not\in C(u_3)$. If $5\not\in C(u_1)$,  let $(uu_1, uu_2,
uu_3, uv)\to (5, 3, 4, 2)$. So assume that $5\in C(u_1)$.

Assume that $4\not\in C(u_1)$. Note that $G$ contains a $(4,
6)_{(u_2, v_2)}$-path and a $(3, i)_{(v_1, v_3)}$-path for any $i\in
\{1, 6\}\backslash C(v_3)$ by ($*_{\textbf{2.1}}$) and
($*_{\textbf{2.2}}$). If $\{1, 6\}\backslash C(v_3)\neq \emptyset$,
then let $(uu_1, uu_2)\to (4, 3)$ and $uv\to \alpha\in \{1,
6\}\backslash C(v_3)$. Otherwise, $C(v_3)= \{5, 3, 1, 6\}$, let
$(uu_1$, $uu_2$, $uv)$ $\to$ $(4, 3, 2)$.

It remains to consider the case $C(u_1)= \{2, 3, 5, 4\}$ and $C(u_3)= \{5, 2, 6, 1\}$
since $\{1, 6\}\subseteq C(u_1)\cup C(u_3)$ by ($*_{\textbf{2.2}}$).
If $2\not\in C(v_3)$,   let $(vv_3, uu_3, uv)\to (2, 3, 5)$.
So assume that $2\in C(v_3)$. If $1, 6\not\in C(v_3)$, then $C(v_1)=
\{3, 2, 1, 6\}$, let $(vv_1, vv_2)\to (4, 2)$. Otherwise,
$\{1, 6\}\cap C(v_3)\neq \emptyset$ and $4\not\in C(v_3)$. If
$4\not\in C(v_1)$, let $(vv_1, vv_2)\to (4, 2)$; otherwise,
$4\in C(v_1)$ and $5\not\in C(v_1)$,  let $(vv_1, vv_2, vv_3,
uv)\to (5, 2, 4, 3)$.

Hence, $4\in C(u_3)\cap C(v_3)$. If $3\not\in C(u_3)$ and
$2\not\in C(v_3)$,  let $(uu_3$, $vv_3$, $uv)$ $\to$ $(3, 2,
5)$. Otherwise, $3\in C(u_3)$ or $2\in C(v_3)$ and assume that
$C(u_3)= \{5, 2, 4, 3\}$ and $C(u_1)= \{3, 2, 1, 6\}$ by
($*_{\textbf{2.2}}$). Let $(uu_1, uu_3, uv)\to (5, 1, 2)$.

\medskip
 \noindent{\bf Case 2.4}\ $C(u_2)\backslash \{1, 4\}= \{5, 6\}$.
\medskip

If $C(v_2)\backslash\{1, 4\}\neq \{5, 6\}$, then we reduce the proof to
 Cases 2.1-2.3. Thus, $C(v_2)\backslash\{1, 4\}$ $= \{5,
6\}$.

\medskip
 \noindent{\bf Case 2.4.1}\   $G$ contains no $(5, 3)_{(u_2,
u_3)}$-path.

\medskip

First, let $uu_2\to 3$. If $G$ contains no $(1, 3)_{(u_2,
v_2)}$-path,  we are done. Otherwise, $G$ contains a $(1,
3)_{(u_2, v_2)}$-path and $1\in C(v_3)$ by ($*_{\textbf{2.2}}$). If
$C(v_1)\backslash \{1, 3\}\neq \{5, 6\}$, then by Cases 1, 2.1-2.3,
we can break this $(1, 3)_{(u, v)}$-cycle and hence break the cycle $B$ so that no new bichromatic cycles are produced.
Thus, assume that $C(v_1)\backslash \{1, 3\}= \{5, 6\}$. Note that
$G$ contains a $(4, 6)_{(u_2, v_2)}$-path and a $(3, 6)_{(v_1,
v_3)}$-path by ($*_{\textbf{2.1}}$) and ($*_{\textbf{2.2}}$). If
$6\not\in C(v_3)$,  let $uv\to 6$. Otherwise, $6\in C(v_3)$. If
$3, 4\not\in C(v_3)$, let $(vv_1, vv_2)= (4, 3)$; otherwise, $\{3,
4\}\cap C(v_3)\neq \emptyset$. If $C(v_3)= \{5, 1, 6, 3\}$,  let
$(vv_2, uv)= (2, 4)$; otherwise, $C(v_3)= \{5, 1, 6, 4\}$,  assuming that
$(uu_2, vv_1, uv)= (4, 2, 3)$, $C(v_3)= \{5, 1, 6, 3\}$,  let $(vv_2, uv)\to (2, 4)$.

\medskip
 \noindent{\bf Case 2.4.2}\  $G$ contains a $(3, 5)_{(u_2,
u_3)}$-path and $(2, 5)_{(v_2, v_3)}$-path.
\medskip

If $2\not\in C(v_1)$, let $(vv_1, uv)\to (2, 3)$; otherwise, $2\in
C(v_1)$. Note that $\{1, 6\}\subseteq C(v_3)$ or $G$ contains a $(3,
i)_{(v_1, v_3)}$-path for any $i\in \{1, 6\}\backslash C(v_3)$ by
($*_{\textbf{2.2}}$).

Assume that $3\not\in C(v_3)$ and then $C(v_3)= \{5, 2, 1, 6\}$ by
($*_{\textbf{2.2}}$). If $4\not\in C(v_1)$,   let $(vv_1,
vv_2)\to (4, 3)$. If $1\not\in C(v_1)$,  let $(vv_1, vv_3,
uv)\to (1, 3, 6)$ and we are done by
($\divideontimes_{\textbf{2}}$). Otherwise, $C(v_1)= \{3, 2, 4, 1\}$
and let $(vv_1, vv_2)\to (6, 3)$.

Assume that $3\in C(v_3)$ and $\{1, 6\}\backslash C(v_3)\neq
\emptyset$. If $5\not\in C(v_1)$,   let $(vv_1, vv_2)\to (5, 3)$,
$vv_3\to \beta_1\in \{1, 6\}\backslash C(v_3)$, $uv\to \{1,
6\}\backslash \beta_1$. Otherwise, $5\in C(v_1)$. Since $\{1,
6\}\subseteq C(v_1)\cup C(v_3)$ by ($*_{\textbf{2.2}}$), we have
$4\not\in C(v_1)\cup C(v_3)$. Let $(vv_1, vv_2)\to (4, 3)$,
$vv_3\to \beta_2\in \{1, 6\}\backslash C(v_3)$, $uv\to \{1,
6\}\backslash \beta_2$.

\medskip

 \noindent {\bf Case 3}\  $|C(u)\cap C(v)|= 4$, say $(uu_1, uu_3)_c= (2, 5)$ and $(vv_1,
vv_3)_c=(2,5)$.

\medskip

If $C(u_2)\backslash$ $\{1, 4\}\neq \{2, 5\}$ or
$C(v_2)\backslash \{1, 4\}\neq \{2, 5\}$, then $|C(u_2)\cap
C(u)|\leq 3$ or $|C(v_2)\cap C(v)|\leq 3$ and we reduce the proof to
 Cases 1 or 2. Thus, $C(u_2)= C(v_2)= \{1, 2, 4, 5\}$. If
$G$ contains neither a $(2, j)_{(u, v)}$-path nor a $(5, j)_{(u,
v)}$-path for some $j\in \{3, 6\}$, then let $uv\to j$. Otherwise,
for any $j\in \{3, 6\}$, $G$ contains an $(i, j)_{(u, v)}$-path for
some $i\in \{2, 5\}$, and $G$ contains a $(2, 3)_{(u, v)}$-path,
$3\in C(u_1)$. If $G$ contains no $(3, 5)_{(v_2, v_3)}$-path,
let $vv_2\to 3$. Otherwise, $G$ contains a $(3, 5)_{(v_2,
v_3)}$-path and a $(3, 5)_{(u_2, u_3)}$-path similarly. It follows
that $3\in C(u_3)$ and $6\in C(u_1)\cap C(u_3)$ similarly. If
$1\not\in C(u_1)$ and $G$ contains no $(1, 5)_{(u_1, u_3)}$-path,
then let $(uu_1, uv)\to (1, 3)$. Hence, $1\in C(u_1)$ or $G$
contains a $(1, 5)_{(u_1, u_3)}$-path. It follows that $\{1, 5\}\cap
C(u_1)\neq \emptyset$, $4\not\in C(u_1)$,  and if $1\not\in C(u_1)$,
then $1\in C(u_3)$. If $4\not\in C(u_3)$,   let $(uu_2, uu_3)\to
(3, 4)$. Otherwise, $C(u_3)= \{5, 3, 6, 4\}$ and $C(u_1)= \{2, 3, 6,
1\}$. Let $(uu_1, uu_3, uv)\to (5, 1, 3)$. \qed

\begin{lemma}\label{4}(\cite{MLSC09})
If $G$ is a graph with $\Delta=4$ that is not 4-regular, then
$a'(G)\le 6$.
 \end{lemma}

\begin{theorem}\label{degree4}
If $G$ is a $4$-regular graph, then $a'(G)\leq 6$.
\end{theorem}

\proof  The proof is proceeded by induction on the vertex number of
$G$. If $|V(G)|=5$, then $G$ is the complete graph $K_5$, and it is
easy to show that $a'(G)\le 6$. Let $G$ be a 4-regular graph with
$|V(G)|\ge 6$. Obviously we may assume that $G$ is 2-connected by
Lemma \ref{4}. If $G$ contains no 3-cycles, then $a'(G)\leq 6$ by
Theorem \ref{triangle-free}. So assume that $G$ contains at least
one 3-cycle. For any graph $H$ with $\Delta(H)=4$ and
$|V(H)|<|V(G)|$, by induction hypothesis or Lemma \ref{4}, $H$
admits  an acyclic edge-6-coloring $c$ using the color set
$C=\{1,2,\ldots,6\}$. Let $v\in V(G)$ be a vertex with neighbors
$v_0, v_1, v_2, v_3$ that lies in at least one 3-cycle.
 To extend $c$ to the original graph $G$,
we split the proof into the following cases.





\medskip
\noindent{\bf Case 1} \ $v_1v_2, v_2v_3, v_3v_0,
v_0v_1\in E(G)$, as shown in Fig.\,3.

\medskip

\begin{figure}[hh]
\begin{center}

\includegraphics[width= 5 in]{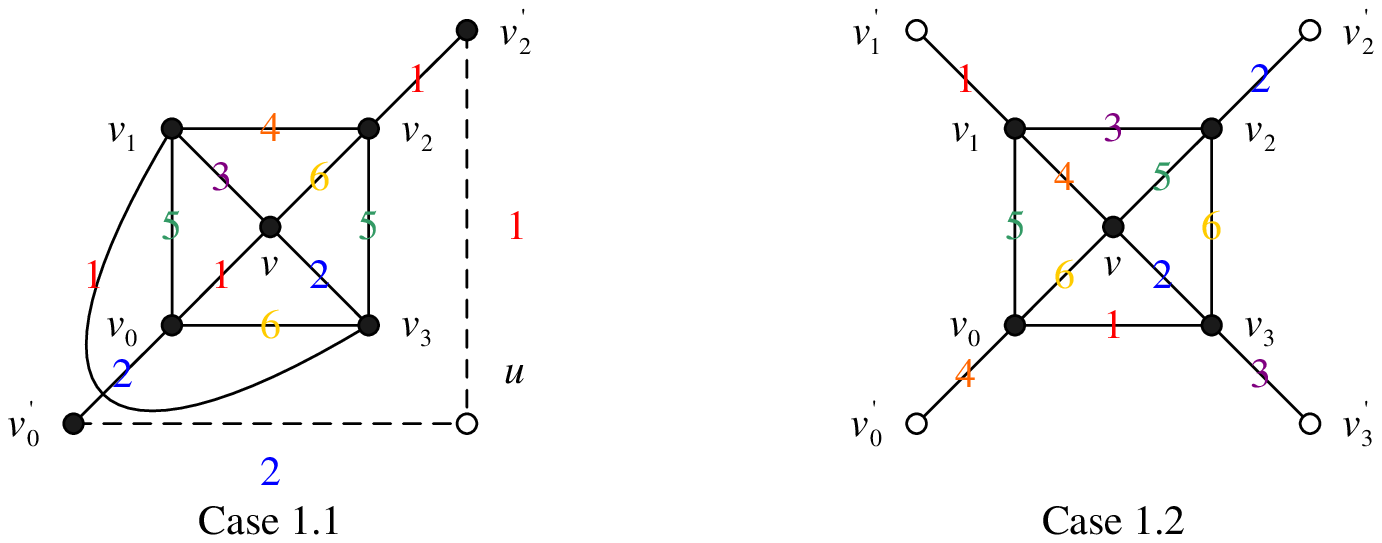}\\

\bigskip

{\rm Fig. 3: The configuration in Case 1.}
\end{center}
\end{figure}

For $i\in \{0,1,2,3\}$, let $v'_i$ be the neighbor of $v_i$
different from $v,v_{i - 1}, v_{i + 1}$, where the indices are taken
modulo $4$. Since $|V(G)|\ge 6$, $G\ne K_5$. We only need to
consider the following subcases by symmetry:

\medskip
 \noindent{\bf Case 1.1}\   $v_1v_3\in E(G)$ and $v_2v_0\not\in
E(G)$.
\medskip

Since $G$ is $2$-connected, $v'_0\neq v'_2$. Let $H= G- \{v, v_1,
v_2, v_3, v_0\} + \{u, uv'_2, uv'_0\}$.  Without loss of generality,
assume that $c(uv'_2)=1$ and $c(uv'_0)=2$.  To extend $c$ to the
graph $G$, we let $\{v_1v_3, v_2v'_2, vv_0\}\to 1$, $\{vv_3,
v_0v'_0\}\to 2$, $\{v_1v_0, v_2v_3\}\to 5$, $\{v_3v_0, vv_2\}\to 6$,
and $(vv_1, v_1v_2)\to (3, 4)$.

\medskip
 \noindent{\bf Case 1.2}\  $v_1v_3,v_2v_0\not\in
E(G)$.
\medskip

If $v'_0,v'_1,v'_2,v'_3$ are identical into a vertex, then
$|V(G)|=6$ and it is easy to give an acyclic edge-6-coloring of $G$.

If $v'_1,v'_2,v'_3$ are identical to a vertex, say $v_5$, and
 $v'_0\neq v_5$, then let $H= G-\{v, v_1, v_2, v_3, v_0\} + \{u, uv'_0, uv_5\}$ with
 $(uv_5, uv'_0)_c= (1, 4)$, $2, 3\not\in C(v_5)$, where $u$ is a new
 vertex added in the graph.

If $v'_1,v'_2$ are are identical to a vertex $v_5$ and $v'_0,v'_3$
are are identical into a vertex $v_6$, we let $H= G- \{v, v_1, v_2,
v_3, v_0\} + \{u, uv_5, uv_6\}$. Assume that $(uv_5, uv_6)_c= (1,
4)$, $2\in C\backslash (C(v_5)\cup \{1, 4\})$, and $3\in C\backslash
(C(v_6)\cup \{1, 4, 2\})$.

If $v'_1,v'_2$ are are identical to a vertex $v_5$, and
  $v'_3,v'_0\ne  v_5$,   let $H= G- \{v, v_1, v_2, v_3, v_0\} + \{u,
uv_5, uv'_3, uv'_0\}$ and assume that $(uv_5, uv'_3, uv'_0)_c= (1,
3, 4)$, $2\in C \backslash (C(v_5)\cup \{3, 4\})$.

 If $v'_1= v'_3= v_5$
and $v'_2= v'_0= v_6$,  let $H= G- \{v,   v_1, v_2, v_3, v_0\} +
\{u, uv_5, uv_6\}$  and assume that $(uv_5, uv_6)_c= (1, 2)$, $3\in
C\backslash (C(v_5)\cup \{1, 2\})$, and $4\in C\backslash
(C(v_6)\cup \{1, 2, 3\})$.

If $v'_1= v'_3= v_5$ and $v'_2\neq v'_0\neq v_5$,  let $H= G- \{v,
v_1, v_2, v_3, v_0\} + \{u, uv_5, uv'_2, uv'_0\}$ and assume that
$(uv_5, uv'_2, uv'_0)_c= (1, 2, 4)$, $3\in C\backslash (C(v_5)\cup
\{1, 2\})$.

If $v'_0,v'_1,v'_2,v'_3$ are all distinct,  let $H= G- \{v, v_1,
v_2, v_3, v_0\} + \{u, uv'_1, uv'_2, uv'_3, uv'_0\}$ and assume that
$(uv'_1, uv'_2, uv'_3, uv'_0)_c= (1, 2, 3, 4)$.

For all cases above, we first let $v_iv'_i\to i$ for $i= 1, 2, 3$ and $v_0v'_0\to 4$.
Then let $\{vv_2, v_1v_0\}\to 5$, $\{vv_0, v_2v_3\}\to 6$, and
$(v_3v_0, vv_3, v_1v_2, vv_1)\to (1, 2, 3, 4)$.

\medskip
\noindent{\bf Case 2}\   $v_0v_1, v_1v_2, v_2v_3\in E(G)$
and $v_3v_0\not\in E(G)$, as shown in  Fig. 4.

\begin{figure}[hh]
\begin{center}

\includegraphics[width= 4.5 in]{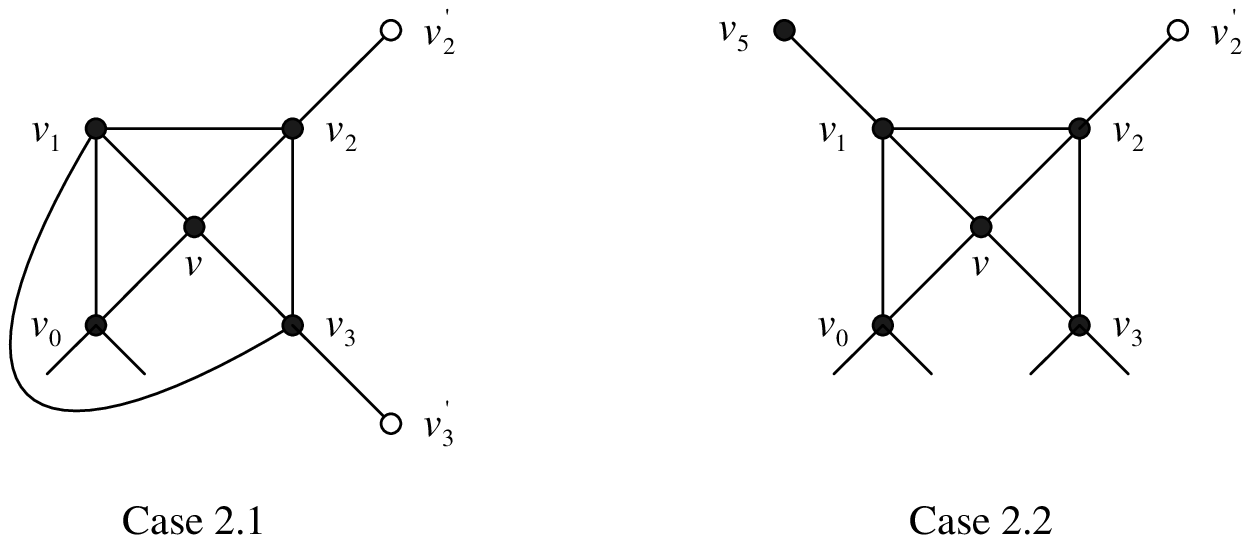}\\

\bigskip

{\rm Fig.\,4: The configuration in Case 2.}
\end{center}
\end{figure}

Let $v'_2$ be the neighbor of $v_2$ different from $v,v_1,v_3$. If
$v_1v_3, v_2v_0\in E(G)$, then $v_1v_0, v_0v_2, v_2v_3, v_3v_1\in
E(G)$ and the proof is reduced to Case 1. Without loss of
generality, assume that $v_2v_0\not\in E(G)$.

\medskip
 \noindent{\bf Case 2.1}\ $v_1v_3\in E(G)$.
\medskip

  Let
$v'_3$ be the neighbor of $v_3$ different from $v,v_1,v_2$.

 If
$v'_2,v'_3$ are identical  to a vertex, say $v_5$, then let $H= G-
\{v, v_1, v_2, v_3\} + \{u, uv_5, uv_0\}$ and assume that $(uv_5,
uv_0)_c= (1, 2)$, $5\in C\backslash (C(v_5)\cup \{1, 2\})$, and
$6\in C\backslash (C(v_0)\cup \{1, 2, 5\})$.

If $v'_3 \neq v'_2$,  let $H= G- \{v, v_1, v_2, v_3\} + \{u, uv_0,
uv'_2, uv'_3\}$ and assume that $(uv_0, uv'_2, uv'_3)_c$ $= (2, 5,
1)$ and $6\in C\backslash (C(v_0)\cup \{1, 2, 5\})$.

For two cases above, we further let $\{v_1v_2, v_3v'_3\}\to 1$, $\{v_2v_3, v_1v_0\}\to 2$,
$\{v_2v_5, vv_3\}\to 5$, $\{vv_0, v_1v_3\}\to 6$, and $(vv_2,
vv_1)\to (3, 4)$.

\medskip
 \noindent{\bf Case 2.2}\  $v_1v_3\not\in E(G)$.
\medskip

 Let $v_5$ be the neighbor of $v_1$ different from $v,v_0,v_2$.
 We consider two subcases as follows:

\medskip
 \noindent{\bf Case 2.2.1}\ $v'_2= v_5$.
\medskip

By Case 1, $v_5v_3, v_5v_0\not\in E(G)$. Let $H= G- \{v, v_1, v_2\}
+ \{v_5v_3, v_5v_0\}$ and assume that $(v_5v_3, v_5v_0)_c= (2, 1)$.
Then, first, let $\{v_1v_5, vv_0\}\to 1$, $\{v_2v_5, vv_3\}\to 2$.
Next, if $\{3, 4, 5, 6\}\backslash (C(v_3)\cup C(v_0))\neq\emptyset$
and assume that $3\not\in C(v_3)\cup C(v_0)$, $4\not\in C(v_0)$,
then let $\{v_1v_0, v_2v_3\}\to 3$ and $(v_1v_2, vv_1, vv_2)\to (4,
5, 6)$. Otherwise, $C(v_0)= \{3, 4\}$ and $C(v_3)= \{5, 6\}$. Then,
let $\{v_1v_0, vv_2\}\to 6$ and $(v_1v_2, v_2v_3, vv_1)\to (3, 4,
5)$. Since $H$ has no bichromatic cycles, it is obvious that $H$ has
no $(1, 2)_{(v_0, v_3)}$-path and thus $G$ contains no $(1,
2)_{(v_0, v)}$-cycle. (See Fig.\,5. We neglect to explain why such
bichromatic cycles cannot exist in the following similar cases.)

\begin{figure}[hh]
\begin{center}

\includegraphics[width= 3.5 in]{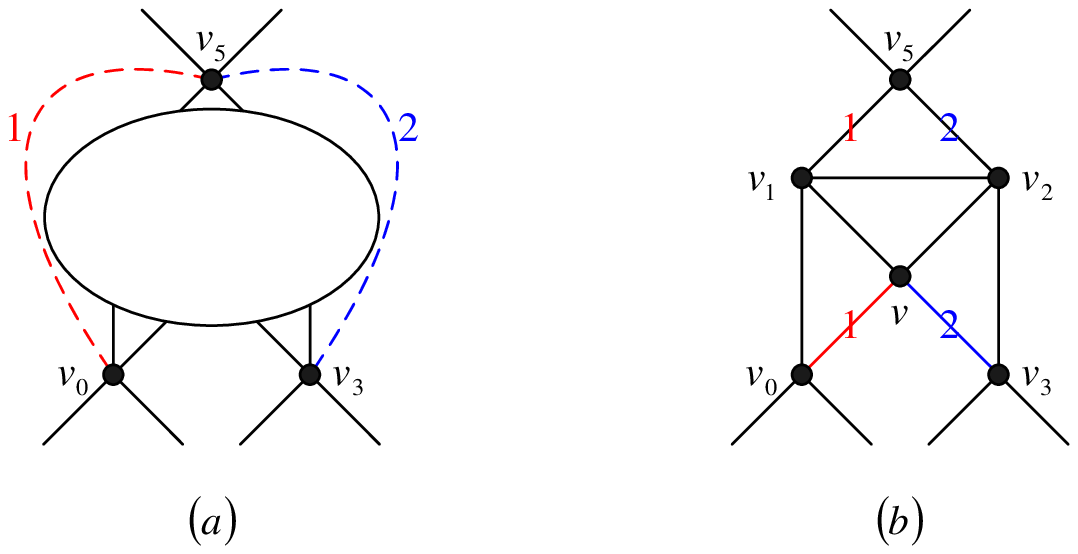}\\

\bigskip

{\rm Fig.\,5: Case 2.2.1.}
\end{center}
\end{figure}

\medskip
 \noindent{\bf Case 2.2.2}\   $v'_2\neq v_5$.
\medskip

Let $H= G- \{v, v_1, v_2\} + \{u, uv_5, uv'_2, uv_3, uv_0, v_3v_0\}$
and assume that $(uv_0, uv_3, v_3v_0)_c$ $=$ $(1, 2, 3)$. Then,
first let $(v_1v_5, v_2v'_2)\to (c(uv_5), c(uv'_2))$.

$\bullet$\ Assume that $(uv_5, uv'_2)= (3, 4)$.

If $\{5, 6\}\backslash C(v_3)\neq \emptyset$ and $5\not\in C(v_3)$,
then let $\{v_1v_0, vv_2\}\to 1$, $\{vv_1, v_2v_3\}\to 2$, and
$(vv_0, vv_3, v_1v_2)\to (3, 5, 6)$. Otherwise, $C(v_3)= \{5, 6\}$.
Then, first let $\{vv_0, v_1v_2\}\to 1$, $(vv_3, v_2v_3)\to (2, 3)$.
And next, if $4\not\in C(v_0)$, then let $(v_1v_0, vv_1, vv_2)\to
(4, 5, 6)$. Otherwise, $4\in C(v_0)$ and assume that $5\not\in
C(v_0)$, then let $(vv_1, v_1v_0, vv_2)\to (4, 5, 6)$.

$\bullet$\ Assume that $(uv_5, uv'_2)= (4, 5)$.

If $4\not\in C(v_0)$, then let $\{vv_0, v_1v_2\}\to 1$, $\{v_1v_0,
vv_3\}\to 3$, and $(v_2v_3, vv_2, vv_1)\to (2, 4, 6)$. Since $H$ has
no bichromatic cycles and $c(v_0v_3)= 3$, it is obvious that $H$
has no $(3, i)_{(v_0, v_3)}$-path for any $i\in C\backslash \{3\}$
and thus $G$ contains no $(3, 6)_{(v_0, v_1)}$-cycle. (See (a), (b)
in Fig.\,6.)

If $4\in C(v_0)$, then let $\{vv_1, v_2v_3\}\to 2$, $(v_1v_0,
vv_0)\to (1, 3)$. If $4\not\in C(v_3)$, then let $(vv_2, vv_3,
v_1v_2)\to (1,4, 6)$. Otherwise, $4\in C(v_3)$ and $5\in C(v_3)\cap
C(v_0)$ similarly. It follows that $C(v_3)= C(v_0)= \{4, 5\}$. Then,
let $(vv_3, v_1v_2, vv_2)\to (1, 3, 6)$. Since $c(v_0v_3)= 3$ in
$H$, it is obvious that $3\not\in C(v_3)$ and thus $G$ contains no
$(3, 2)_{(v_0, v)}$-cycle. (See (a), (c) in Fig.\,6.)

\begin{figure}[hh]
\begin{center}

\includegraphics[width= 5.5 in]{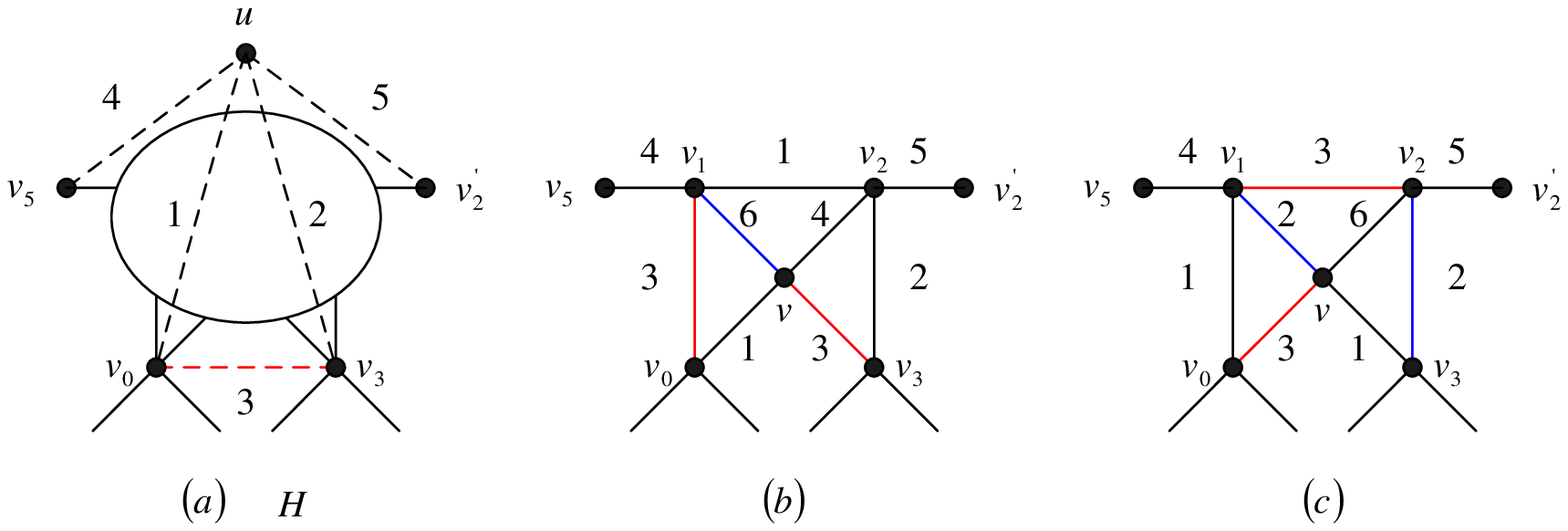}\\

\bigskip

{\rm Fig.\,6: Case 2.2.2.}
\end{center}
\end{figure}

\medskip

\noindent{\bf Case 3}\  $v_1v_2, v_2v_3\in E(G)$ and
$v_1v_0, v_3v_0\not\in E(G)$, as shown in Fig.\,7.

\medskip

\begin{figure}[hh]
\begin{center}

\includegraphics[width= 5.5 in]{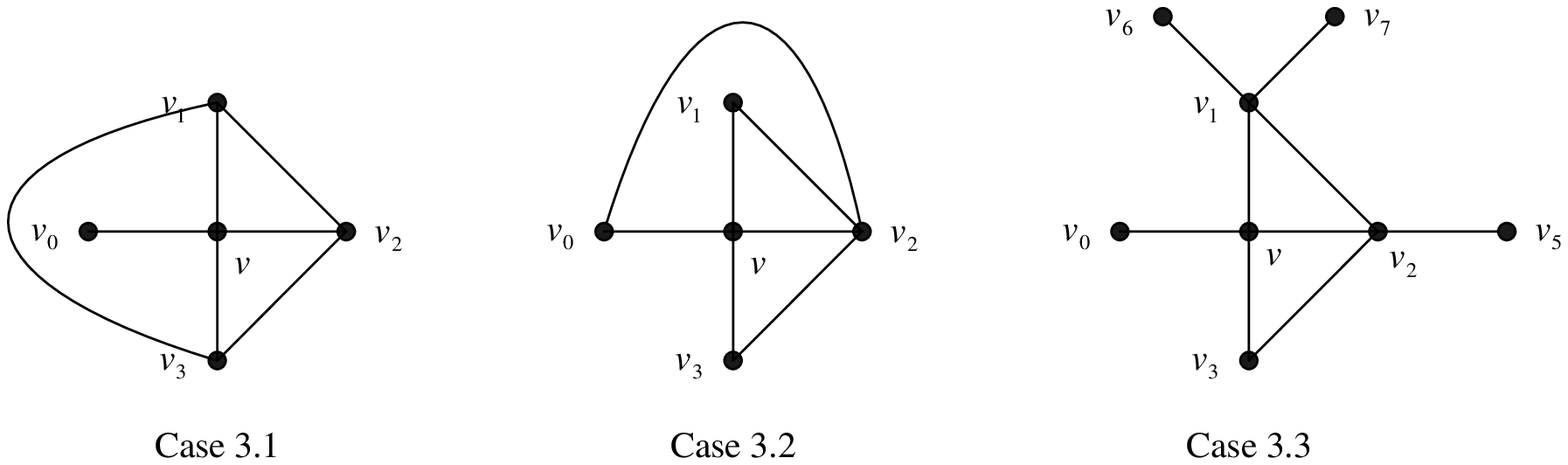}\\

\bigskip

{\rm Fig.\,7: The configuration in Case 3.}
\end{center}
\end{figure}

If $v_1v_3, v_2v_0\in E(G)$, then $v_1v_3, v_3v_2, v_2v_0\in E(G)$
and the proof can be reduced to  Case 2.
Otherwise, $v_1v_3\not\in E(G)$ or $v_2v_0\not\in E(G)$.

\medskip
\noindent{\bf Case 3.1}\  $v_1v_3\in E(G)$ and $v_2v_0\not\in
E(G)$.
\medskip

Let $v_5$ be the neighbor of $v_1$ different from $v,v_2,v_3$. Let
$v_6$ be the neighbor of $v_2$ different from $v,v_1,v_3$. Let $v_7$
be the neighbor of $v_3$ different from $v,v_1,v_2$. By Cases 1 and
2, we can assume that $v_0, v_5, v_6, v_7$ are pairwise distinct.

Let $H= G- \{v, v_1, v_2, v_3\} + \{u, uv_0, uv_5, uv_6, uv_7\}$ and
assume that $(uv_5, uv_6, uv_7, uv_0)_c$ $=$ $(1, 2, 3, 4)$. Then,
let $\{v_1v_5, v_2v_3\}\to 1$, $\{vv_3, v_2v_6\}\to 2$, $\{vv_1,
v_3v_7\}$ $\to 3$, $\{vv_2, v_1v_3\}$ $\to$ $5$, and $(vv_0,
v_1v_2)\to (4, 6)$.

\medskip
 \noindent{\bf Case 3.2}\  $v_1v_3\not\in E(G)$ and $v_2v_0\in
E(G)$.

\medskip

Let $H= G- \{v, v_2\} + \{v_1v_0, v_0v_3, v_3v_1\}$ and assume that
$(v_1v_0, v_0v_3, v_3v_1)_c= (1, 2, 3)$, $4\in C\backslash
(C(v_1)\cup \{1, 2, 3\})$. Then, first let $\{vv_3, v_1v_2\}\to 3$,
$(vv_0, vv_1)\to (1, 4)$. Next, if $\{5, 6\}\backslash C(v_3)\neq
\emptyset$ and $5\not\in C(v_3)$, then let $(v_2v_0, v_2v_3,
vv_2)\to (2, 5, 6)$; otherwise, $C(v_3)= \{5, 6\}$ and let
$v_2v_3\to 1$. If $\{5, 6\}\backslash C(v_0)\neq \emptyset$, then
let $v_2v_0\to 2$ and $vv_2\to\beta\in \{5, 6\}\backslash C(v_0)$;
otherwise, $C(v_0)= \{5, 6\}$ and let $(vv_2, v_2v_0)\to (2, 4)$.

\medskip
 \noindent{\bf Case 3.3}\   $v_1v_3\not\in E(G)$ and $v_2v_0\not\in
E(G)$.
\medskip

Let $v_5$ be the neighbor of $v_2$ different from $v,v_1,v_3$.  Let
$v_6, v_7$ be the other two neighbors of $v_1$ different from
$v,v_2$.  Note that
 $v_0, v_5, v_6, v_7$ are pairwise distinct by the previous argument. Let
$H= G- \{v, v_2\} + \{v_1v_0, v_1v_5\}$ and assume that $(v_1v_5$,
$v_1v_0$, $v_1v_6$, $v_1v_7)_c$ $=$ $(1, 2, 3, 4)$. Then, first let
$(v_2v_5, vv_0)\to (1, 2)$.

\medskip
 \noindent{\bf Case 3.3.1}\ $C(v_3)= \{1, 2\}$.
\medskip

Assume that $6\not\in C(v_5)$. Then, first let $(v_1v_2, vv_2,
v_2v_3)\to (2, 4, 6)$. Next, if $5\not\in C(v_0)$, then let $(vv_1,
vv_3)\to (1, 5)$; if $3\not\in C(v_0)$, then let $(vv_1, vv_3)\to
(5, 3)$; otherwise, $C(v_0)= \{2, 5, 3, 4\}$ and let $(vv_1, vv_3,
vv_0)\to (1, 5, 6)$.

Assume that $\{5, 6\}\subseteq C(v_5)$ and $\{5, 6\}\subseteq
C(v_0)$ similarly. If $3, 4\not\in C(v_0)\cup C(v_5)$, then let
$\{vv_1, v_2v_3\}\to 5$ and $(v_1v_2, v_2v_5, vv_3, vv_2)\to (1, 3,
4, 6)$; otherwise, $C(v_0)= \{2, 5, 6, 4\}$. Then, let $(vv_1, vv_2,
vv_0, vv_3, v_1v_2)\to (1, 2, 3, 5, 6)$ and $v_2v_3\to\gamma \in
\{3, 4\}\backslash C(v_5)$.

\medskip
 \noindent{\bf Case 3.3.2}\  $C(v_3)= \{1, 3\}$.
\medskip

If $6\not\in C(v_5)$, then let $\{vv_1, v_2v_3\}\to 6$ and $(v_1v_2,
vv_3, vv_2)$ $\to$ $(2, 4, 5)$. If $G$ contains no $(1, 4)_{(v_3,
v_5)}$-path, then let $(vv_1, vv_2, v_2v_3, vv_3, v_1v_2)\to (1, 3,
4, 5, 6)$. Otherwise, $G$ contains a $(1, 4)_{(v_3, v_5)}$-path and
hence  $C(v_5)= \{1, 5, 6, 4\}$. We let $(vv_1, v_2v_3, vv_3,
v_1v_2, vv_2)$ $\to$ $(1, 2, 4, 5, 6)$.

\medskip
 \noindent{\bf Case 3.3.3}\  $C(v_3)= \{1, 5\}$.
\medskip

Let $vv_1\to 6$. If $6\not\in C(v_5)$, we further let $(v_1v_2,
vv_3, vv_2$, $v_2v_3)$ $\to$ $(2, 3, 5, 6)$; If $3\not\in C(v_5)$,
  let $(v_1v_2, v_2v_3, vv_3, vv_2)\to (2, 3, 4, 5)$. Otherwise,
$C(v_5)= \{1, 6, 3, 4\}$,  let $(v_2v_3, vv_3, vv_2, v_1v_2)\to (2,
3, 4, 5)$.

\medskip
 \noindent{\bf Case 3.3.4}\ $C(v_3)\notin \{  \{1, 2\},  \{1, 3\},   \{1, 5\}\}$.

\medskip

If $C(v_3)= \{5, 6\}$,  let $(vv_1, v_2v_3, vv_3, vv_2, v_1v_2)\to
(1, 2, 3, 4, 5)$. If $C(v_3)$ $=$ $\{3, 5\}$,   let $\{vv_1,
v_2v_3\}\to 6$ and $(vv_3, vv_2, v_1v_2)\to (1, 4, 5)$. If $C(v_3)=
\{3, 4\}$,   let $\{vv_1, v_2v_3\}\to 5$ and $(vv_3, v_1v_2,
vv_2)\to (1, 2, 6)$.

\medskip
\noindent{\bf Case 4}\   $v_1v_2, v_3v_0\in E(G)$ and $v_0v_1,
v_2v_3, v_1v_3, v_2v_0\not\in E(G)$.
\medskip

Let $V_i= \{v_{i1}, v_{i2}\}$ be the set of the other neighbors of
$v_i$ for $i\in \{0, 1, 2, 3\}$. By Cases 1 to 3, $V_1\cap V_2=
\emptyset$ and  $V_3\cap V_0= \emptyset$. Let $H= G- \{v, v_1, v_2,
v_3, v_0\} + \{u, w, uv_{11}, uv_{12}, uv_{21}, uv_{22}, wv_{31},
wv_{32}, wv_{01}, wv_{02}\}$, where $u$ and $w$ are new vertices
added. Assume that $(uv_{11}$, $uv_{12}$, $uv_{21}$, $uv_{22})_c$
$=$ $(1$, $2$, $3$, $4)$, $\{c(wv_{01}), c(wv_{02})\}$ $=$ $\{a,
b\}$, $\{c(wv_{31}), c(wv_{32})\}= \{c, d\}$, where $\{a, b\}\cap
\{c, d\}= \emptyset$, see Fig.\,8.  Then, first let $v_iv_{ij}\to
c(uv_{ij})$ for $i\in\{ 1, 2\}$, and $vv_{ij}\to c(wv_{ij})$ for
$i\in \{3, 4\}$, $j\in \{1, 2\}$. Next, by symmetry, we need to
consider the following.

\begin{figure}[hh]
\begin{center}

\includegraphics[width= 5 in]{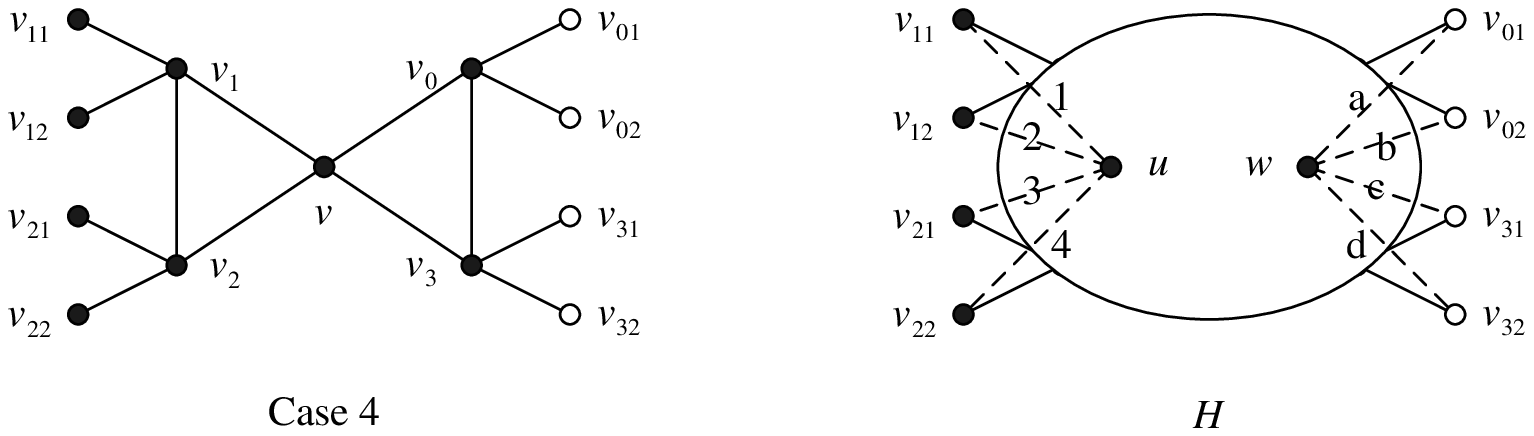}\\

\bigskip

{\rm Fig.\,8: The configuration in Case 4.}
\end{center}
\end{figure}
\medskip
\noindent{\bf Case 4.1}\   $\{a, b\}= \{1, 3\}$ and $\{c, d\}= \{2,
4\}$.
\medskip

If $5\not\in C(v_{11})$,  let $\{v_1v_2, vv_0\}\to 5$ and $(vv_2,
vv_3, vv_1, v_3v_0)\to (1, 3, 4, 6)$. If $4\not\in C(v_{11})$,  let
$\{v_1v_2, vv_0\}\to 5$ and $(vv_3, vv_2, vv_1, v_3v_0)\to (1, 2, 4,
6)$. Otherwise, $C(v_{11})= \{1, 4, 5, 6\}$ and $C(v_{12})= \{2, 3,
5, 6\}$. Let $\{v_1v_2, vv_3\}\to 1$, $\{vv_1, v_3v_0\}\to 5$, and
$(vv_2, v_1v_{11}, vv_0)\to (2, 3, 6)$.

\medskip
 \noindent{\bf Case 4.2}\   $\{a, b\}= \{1, 2\}$ and $\{c, d\}= \{3,
4\}$.
\medskip

If $G$ contains no $(i, j)_{(v_2, v_0)}$-path for some $i\in \{1,
2\}$ and $j\in \{3, 4\}$, then let $\{v_1v_2, vv_3\}\to 5$, $\{vv_1,
v_3v_0\}\to 6$, and $(vv_2, vv_0)\to (i, j)$. Otherwise, $G$
contains a $(i, j)_{(v_2, v_0)}$-path for any $i\in \{1, 2\}$ and
$j\in \{3, 4\}$ and $\{1, 2\}\subseteq C(v_{21})\cap C(v_{22})$. It
follows that $\{5, 6\}\backslash C(v_{21})\neq \emptyset$ and assume
that $5\not\in C(v_{21})$. First, let $\{vv_2, v_3v_0\}\to 5$,
$\{v_1v_2, vv_3\}\to 6$, and $(vv_0, vv_1)\to (3, 4)$. Next, if
$6\not\in C(v_{22})$, then we are done; otherwise, $C(v_{22})= \{4,
1, 2, 6\}$ and let $\{vv_1, v_2v_{22}\}\to 5$, $vv_2\to 4$.

\medskip
 \noindent{\bf Case 4.3}\ All other cases.
\medskip

If $\{a, b\}= \{1, 2\}$ and $\{c, d\}= \{5, 6\}$, let $\{vv_1,
v_3v_0\}$ $\to$ $3$ and $(vv_2$, $vv_3$, $v_1v_2$, $vv_0)$ $\to$
$(1$, $4$, $5$, $6)$. If $\{a, b\}= \{1, 3\}$ and $\{c, d\}= \{5,
6\}$, let $(vv_2$, $vv_0$, $vv_3$, $v_3v_0, v_1v_2, vv_1)$ $\to$
$(1, 2, 3, 4, 5, 6)$. If $\{a, b\}= \{1, 5\}$ and $\{c, d\}= \{2,
6\}$, let $(vv_2$, $vv_0$, $v_3v_0$, $vv_3, v_1v_2, vv_1)$ $\to$
$(1, 2, 3, 4, 5, 6)$. If $\{a, b\}= \{1, 5\}$ and $\{c, d\}= \{3,
6\}$, let $\{v_1v_2, vv_3\}$ $\to$ $5$ and $(vv_2, vv_0, v_3v_0,
vv_1)\to (1, 2, 4, 6)$. If $\{a, b\}= \{1, 2\}$ and $\{c, d\}= \{3,
5\}$, let $\{v_1v_2, v_3v_0\}\to 6$ and $(vv_2, vv_1, vv_3, vv_0)\to
(1, 3, 4, 5)$. If $\{a, b\}= \{1, 3\}$ and $\{c, d\}= \{2, 5\}$, let
$\{v_1v_2, v_3v_0\}\to 6$ and $(vv_2, vv_0, vv_3, vv_1)\to (1, 2, 4,
5)$.

\medskip
\noindent{\bf Case 5}\   $v_1v_2\in E(G)$ and $v_2v_3, v_3v_0,
v_0v_1, v_1v_3, v_2v_0\not\in E(G)$.
\medskip

Let $V_1,V_2$ be defined similarly as in Case 4.  By Cases 1 to 4  ,
$V_1\cap V_2= \emptyset$. Let $H= G- \{v, v_1, v_2\} + \{u, uv_{11},
uv_{12}, uv_{21}, uv_{22}, v_3v_0\}$ and assume $(uv_{11}, uv_{12},
uv_{21}, uv_{22})_c= (1, 2, 3, 4)$, as shown in Fig.\,9.  Then,
first let $v_iv_{ij}\to c(uv_{ij})$ for $i, j\in\{ 1, 2\}$. Next,
since $c(v_3v_0)= 1$ or $c(v_3v_0)= 5$, we need to consider the
following.

\begin{figure}[hh]
\begin{center}

\includegraphics[width= 5 in]{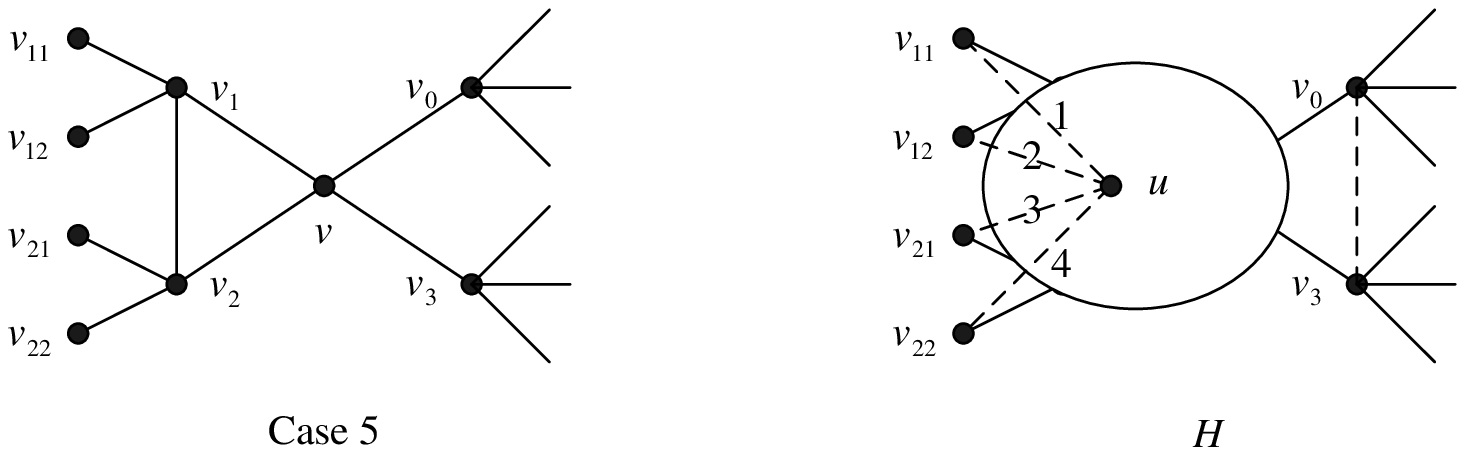}\\

\bigskip

{\rm Fig.\,9:  The configuration in Case 5.}
\end{center}
\end{figure}

\medskip
 \noindent{\bf Case 5.1}\   $c(v_3v_0)= 1$.
\medskip

\noindent {\bf Case 5.1.1}\ $2\not\in C(v_3)$.
\medskip

 If $G$ contains neither  $(2, 3)_{(v_1, v_3)}$-path nor   $(1,
3)_{(v_1, v_0)}$-path,   let $(vv_0, vv_3, vv_1, v_1v_2, vv_2)$ $\to
(1, 2, 3, 5, 6)$. Otherwise, $G$ contains either a $(2, 3)_{(v_1,
v_3)}$-path or a $(1, 3)_{(v_1, v_0)}$-path. If $3\in
C(v_3)\backslash C(v_0)$, then let $(vv_2, vv_3, vv_0)\to (1, 2, 3)$
and $vv_1\to x_1\in \{4, 5, 6\}\backslash C(v_3)$, $v_1v_2\to y_1\in
\{5, 6\}\backslash \{x_1\}$. If $3\in C(v_0)\backslash C(v_3)$, then
let $(vv_0, vv_2, vv_3)\to (1, 2, 3)$ and $vv_1\to x_2\in \{4, 5,
6\}\backslash C(v_0)$, $v_1v_2\to y_2\in \{5, 6\}\backslash
\{x_2\}$. Otherwise, $3\in C(v_3)\cap C(v_0)$, $4\in C(v_3)\cap
C(v_0)$ similarly and then $5\not\in C(v_3)$. First, let $\{v_1v_2,
vv_3\}\to 5$ and $(vv_0, vv_2, vv_1)\to (1, 2, 6)$. Next, if $G$
contains no $(1, 6)_{(v_1, v_0)}$-path, then we are done; otherwise,
$G$ contains a $(1, 6)_{(v_1, v_0)}$-path, which cannot pass through
$v_3$ and $C(v_0)= \{3, 4, 6\}$. We switch the colors of $vv_3$ and
$vv_0$. Otherwise, assume that $2\in C(v_3)$ and $2\in C(v_0)$
similarly.

\medskip
\noindent {\bf Case 5.1.2}\ $3\not\in C(v_3)$.
\medskip

 First, assume that $G$ contains no $(1, 4)_{(v_1, v_0)}$-path and
let $(vv_0, vv_3, vv_1)\to (1, 3, 4)$. If $G$ contains no $(3,
i)_{(v_2, v_3)}$-path for some $i\in \{2, 5, 6\}$, then let $vv_2\to
i$ and   $v_1v_2\to \alpha\in \{5, 6\}\backslash \{i\}$. Otherwise,
$G$ contains a $(3, i)_{(v_2, v_3)}$-path for any $i\in \{2, 5, 6\}$
and $C(v_{21})= \{3, 2, 5, 6\}$, $C(v_3)= \{2, 5, 6\}$. If $\{5,
6\}\backslash C(v_0)\neq \emptyset$ and $5\not\in C(v_0)$, then let
$(v_2v_{21}, v_1v_2, vv_2)\to (1, 3, 5)$; otherwise, $\{5,
6\}\subseteq C(v_0)$, $C(v_0)= \{1, 2, 5, 6\}$ and let $(vv_2, vv_0,
v_1v_2, vv_1)\to (1, 4, 5, 6)$. Next, assume that $G$ contains a
$(1, 4)_{(v_1, v_0)}$-path and $4\in C(v_0)$, $5\not\in C(v_0)$.
Then, let $(vv_3, vv_2, vv_1, vv_0, v_1v_2)\to (1, 2, 4, 5, 6)$.
Otherwise, $3\in C(v_3)$ and $3\in C(v_0)$.

\medskip
\noindent {\bf Case 5.1.3}\ $C(v_3)= C(v_0)= \{2, 3, 4\}$.
\medskip

 First, let $\{v_1v_2, vv_0\}\to 5$ and $(vv_3, vv_2, vv_1)\to (1, 2,
6)$. Next, if $G$ contains no $(2, 5)_{(v_1, v_0)}$-path, then we
are done; otherwise, $G$ contains a $(2, 5)_{(v_1, v_0)}$-path, we
 switch the colors of $vv_3$ and $vv_0$.

\medskip
 \noindent{\bf Case 5.2}\   $c(v_3v_0)= 5$.
\medskip

By symmetry, we may assume that $1\not\in C(v_3)$. If $G$ contains
no $(1, i)_{(v_1, v_3)}$-path for some $i\in \{3, 4\}$, then let
$(vv_3, vv_2, vv_1, vv_0, v_1v_2)\to (1, 2, i, 5, 6)$. Otherwise,
$G$ contains a $(1, i)_{(v_1, v_3)}$-path for any $i\in \{3, 4\}$
and $\{3, 4\}\subseteq C(v_{11})\cap C(v_3)$. If $1\not\in C(v_0)$,
then let $(vv_0, vv_2, vv_1, vv_3, v_1v_2)\to (1, 2, 3, 5, 6)$.
Otherwise, $1\in C(v_0)$.

\medskip
 \noindent{\bf Case 5.2.1}\   $G$ contains no $(2, 5)_{(v_1,
v_0)}$-path.
\medskip

If $G$ contains no $(1, 6)_{(v_1, v_3)}$-path, then let $\{v_1v_2,
vv_0\}\to 5$ and $(vv_3, vv_2, vv_1)\to (1, 2, 6)$. Otherwise, $G$
contains a $(1, 6)_{(v_1, v_3)}$-path and $6\in C(v_{11})\cap
C(v_3)$ and $C(v_{11})= \{1, 3, 4, 6\}$, $C(v_3)= \{3, 4, 6\}$. If
$G$ contains no $(2, i)_{(v_1, v_3)}$-path for some $i\in \{3, 4,
6\}$, then let $(vv_2, vv_3, vv_4, vv_1)\to (1, 2, 5, i)$ and
$v_1v_2\to x_4\in \{5, 6\}\backslash \{i\}$. Otherwise, $G$ contains
a $(2, i)_{(v_1, v_3)}$-path for any $i\in \{3, 4, 6\}$ and
$C(v_{12})= \{2, 3, 4, 6\}$. If $2\not\in C(v_0)$, then let $(vv_2,
vv_0, vv_1, vv_3$, $v_1v_2)$ $\to$ $(1, 2, 3, 5, 6)$. If $6\not\in
C(v_0)$, then let $\{v_1v_2, vv_0\}\to 6$ and $(vv_3, vv_2, vv_1)\to
(1, 2, 5)$. Otherwise, $C(v_0)= \{1, 2, 6\}$ and $4\not\in C(v_0)$.
Then, let $\{v_1v_2, vv_3\}\to 1$, $\{vv_0, v_1v_{11}\}\to 5$, and
$(vv_1, vv_2)\to(4, 6)$.

\medskip
 \noindent{\bf Case 5.2.2}\   $G$ contains a $(2, 5)_{(v_1,
v_0)}$-path, and $5\in C(v_{12})$, $2\in C(v_0)$.
\medskip

Then, $\{3, 4\}\backslash C(v_0)\neq \emptyset$ and assume that
$3\not\in C(v_0)$. If $5\not\in C(v_{11})$, then let $\{vv_0,
v_1v_{11}\}\to 5$ and $(vv_3, vv_2, vv_1, v_1v_2)\to (1, 2, 3, 6)$.
Otherwise, $C(v_{11})= \{1, 3, 4, 5\}$. Then, let $\{v_1v_2,
vv_3\}\to 1$ and $(vv_0, vv_1, vv_2, v_1v_{11})\to (3, 4, 5, 6)$.
\qed


\end{document}